\documentclass[twoside,11pt]{article}
\usepackage{comment}

\usepackage{blindtext}
\usepackage{amsthm}
\usepackage{amssymb}
\usepackage{amsmath}
\usepackage{dsfont}
\usepackage{bbm}
\usepackage{graphicx}
\usepackage{ragged2e}
\usepackage{geometry}
\usepackage[font=scriptsize]{caption}
\usepackage{algorithm}
\usepackage{algorithmic}
\usepackage{chngcntr}
\usepackage{mathtools}
\usepackage[english]{babel}
\usepackage{caption}
\usepackage{float}
\usepackage[round]{natbib}
\newcounter{mycount}
\counterwithin{algorithm}{mycount}
\refstepcounter{mycount}

\newtheorem{remark}{Remark}
\newtheorem{corollary}{Corollary}
\newtheorem{definition}{Definition}
\newtheorem{theorem}{Theorem}
\newtheorem{lemma}{Lemma}
\newtheorem{assumption}{Assumption}

\usepackage[mathscr]{euscript}
\usepackage[colorlinks = true,
            linkcolor = blue,
            urlcolor  = blue,
            citecolor = blue,
            anchorcolor = blue]{hyperref}
\usepackage{soul}
\usepackage{subcaption}
\usepackage[abbrvbib,preprint]{jmlr2e}
\firstpageno{1}

\ShortHeadings{Testing For LLM Response Differences}{A. Acharrya, C.E. Priebe, H.S. Helm}

\firstpageno{1}

\begin{document}


\title{Testing for LLM response differences: \\
the case of a composite null consisting of \\semantically irrelevant query perturbations
}

\author{\name Aranyak Acharyya \email aachary6@jhu.edu \\
       \addr Mathematical Institute for Data Science \\
       Johns Hopkins University\\
       Baltimore, MD 21218, USA
\AND
       \name Carey E. Priebe \email cep@jhu.edu \\
       \addr Department of Applied Mathematics and Statistics \\
       Johns Hopkins University\\
       Baltimore, MD 21218, USA
\AND
       \name Hayden S. Helm \email hayden@helivan.io \\
       \addr Helivan\\
       San Francisco, CA 94123, USA
}


\maketitle

\begin{abstract}
Given an input query, generative models such as large language models 
 produce 
 a random response drawn from a response distribution.
Given two input queries, it is natural to ask if their response distributions are the same.
While traditional statistical hypothesis testing is designed to address this question, the response distribution induced by an input query is often sensitive to semantically irrelevant perturbations to the query, so much so that a traditional test of equality  might indicate that two semantically equivalent queries induce statistically different response distributions.
As a result, the outcome of the statistical test may not align with the user's requirements.
In this paper, we address this misalignment by incorporating into the testing procedure consideration of a collection of semantically similar queries.
In our setting, the mapping from the collection of user-defined semantically similar queries to the corresponding collection of response distributions is not known \textit{a priori} and must be estimated, with a fixed budget.
Although the problem we address is quite general, we focus our analysis on the setting where the responses are binary, show that the proposed test is asymptotically valid and consistent, and discuss important practical considerations with respect to power and computation. 
\end{abstract}


\section{Introduction}
\label{Sec:Introduction}

Our analysis is motivated by a simple observation when working with generative models:
a small change to a query typically changes the response distribution. 
\begin{figure}[t]
    \centering
    \includegraphics[width=\linewidth]{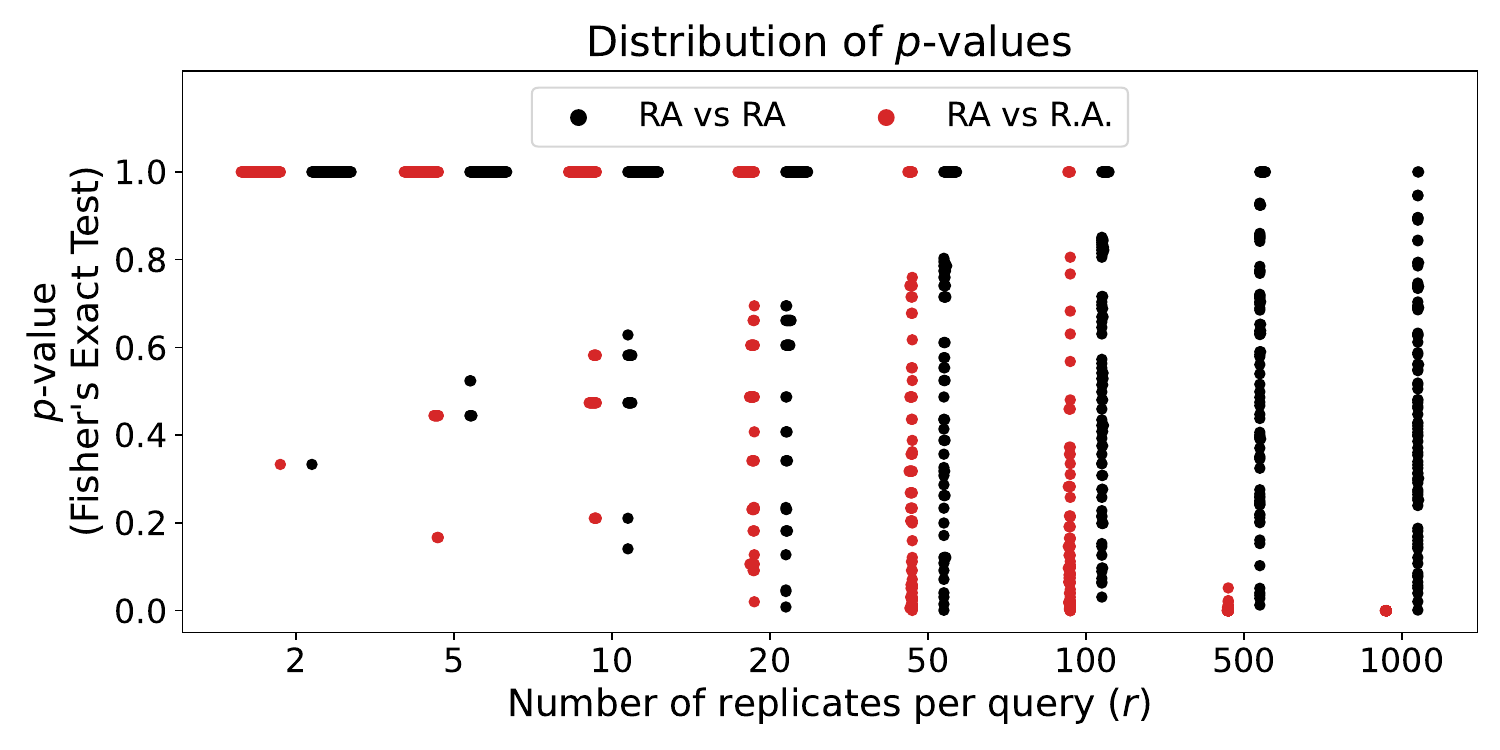}
    \caption{
Empirical distributions for $p$-values  
when testing for equality of the binary response distributions
for
$q_1 =$ ``RA Fisher was a statistician. Was he great?''
against itself (control = black)
or against the semantically irrelevant perturbation
$q_2 =$ ``R.A.\ Fisher was a statistician. Was he a great man?'' (condition = red).
The figure presents $p$-values, for various values of $r$, from Fisher's exact test for 100 Monte Carlo experiments
based on independent samples obtained by repeteadly prompting 
\texttt{LLaMA-3-8B-Instruct}
$r$ times. 
When $r$ is large,
the distribution of $p$-values when introducing a semantically irrelevant change to the query
deviates dramatically from the distribution of $p$-values under the control condition.
With $r=168700$ independent samples for each query,
$\hat{p}_1 \approx 0.870$ and $\hat{p}_2 \approx 0.948$
and $p$-value $\approx 0$;  
while $q_1 \approx q_2$, $p_1 \neq p_2$.
While the user may believe they are under the null in both settings, the sensitivity of response distribution to semantically irrelevant query perturbations produces unwanted rejections (from the perspective of the user) when $r$ is large.
    }
    \label{fig:motivation}
\end{figure}
As an example, consider two nearly identical queries:
$q_1 =$ ``RA Fisher was a statistician. Was he great?''
and
$q_2 =$ ``R.A.\ Fisher was a statistician. Was he great?''
To the majority of English speakers, the two queries are the same --
``RA Fisher'' (without dots) vs.\ ``R.A.\ Fisher'' (with dots) is a semantically irrelevant distinction.
The impact on the 
response distribution, on the other hand, is significant. 
To wit:
letting $p_j$ be the probability that a generative model outputs ``Yes'' in response to $q_j$ for $j=1,2$,
we find that a classical two-sample Neyman-Pearson test for equality of two Bernoulli parameters
$$H_0: p_1 = p_2 \text{~ vs. ~} H_A: p_1 \neq p_2$$ 
yields $p$-value 
close to $0$ given enough samples from the response distributions.
See Figure \ref{fig:motivation} for details.

We view rejections of $H_0$ when $q_1$ and $q_2$ differ by only a semantically irrelevant distinction
to represent {\em statistically} significant findings that are not {\em operationally} significant to the user.
Alas, such rejections are not controlled in the classical Neyman-Pearson testing framework when using a simple null and simple alternative.
Consider the natural composite extension
\begin{align*}
    H_0: p \in \mathcal{P}_0 \text{~ vs. ~} H_A: p \not\in \mathcal{P}_0
\end{align*}
where $\mathcal{P}_0$ is a set of unknown null probabilities induced by a set $\mathcal{Q}_0$ of semantically irrelevant perturbations of a base query $q_0$.
When $ \mathcal{P}_{0} $ is known we can apply composite extensions of our preferred testing procedure. 
In our setting, however, 
the map from $\mathcal{Q}$ to $\mathcal{P}$ is not known, and hence  $\mathcal{P}_{0} $ must be estimated by repeatedly sampling responses from the model for queries in $ \mathcal{Q}_{0}$.
The goal of the current manuscript is to develop a statistical test that
takes into account a set of user-defined semantically similar queries and properly controls
the Type-I error while providing desirable power;
that is, we provide an asymptotically valid and consistent test
for generative model response differences in the case of a composite null consisting of response distributions induced by semantically irrelevant perturbations.

\subsection{Problem Statement}
\label{Subsec:Problem Statement}
For our purposes, a generative model $ f $ is a random mapping from an input space $ \mathcal{Q} $ to an output space $ \mathcal{X} $. 
In particular, given an input (or ``query") $ q \in \mathcal{Q}$ , the random response $ f(q) $ is sampled from a distribution $ F_{f(q)} $ on the set of possible responses.
Repeatedly querying the same model $ r $ times with the same query $ q $ yields $ i.i.d. $ samples $ f(q)_{1}, \hdots, f(q)_{r} $ from $ F_{f(q)} $. 
We let $ g: \mathcal{X} \to \mathbb{R}^{s} $ denote an embedding function that maps from the output space to $ s $-dimensional Euclidean space. 
The embedded response $ g(f(q)) $ is a random vector in $ \mathbb{R}^{s} $ and the replicates $ g(f(q)_1), \dots, g(f(q)_r) $ are $ i.i.d. $ samples from $ F_{g(f(q))} $. 
Due to practical considerations, our analysis is focused on the embedded responses as opposed to distributions on token-strings, and we refer to $ F_{g(f(q))} $ as $ F_{q} $ for notational convenience.

Of primary interest is determining if the two response distributions $ F_{q_{0}} $ and $ F_{q'}$ induced by $ q_{0} $ and $ q' $, respectively, are the same. 
That is,
\begin{equation}
\label{eq:primary-interest}
    H_{0}: F_{q'} = F_{q_0} \quad\quad\quad vs \quad\quad\quad H_{A}: F_{q'} \neq F_{q_0}.
\end{equation}
Given samples $ f(q_{0})_{1}, \hdots, f(q_{0})_{r} $ (respectively, $ f(q')_{1}, \hdots, f(q')_{r} $ ) we can obtain an estimate of $F_{q_0}$, denoted by $\hat{F}_{q_0}$ (respectively, an estimate of $F_{q'}$, denoted by $\hat{F}_{q'}$)
and apply a standard statistical hypothesis test. 
However, as demonstrated by our motivating example above,   standard hypothesis tests in this context may lead to rejections of $ H_0 $ that are not desirable to the user.

To address these operationally 
undesirable
rejections, we define a user-specified set of queries semantically similar to $ q_{0} $, $ \mathcal{Q}_{0} \subseteq \mathcal{Q} $, and modify Eq.\ \eqref{eq:primary-interest}.
Each element  $ q_{i} \in \mathcal{Q}_{0} $ is such that the 
\textit{user} 
expects for an (asymptotically) valid test to have approximately size $ \alpha $ when testing for equality of $ F_{q_{0}} $ and $ F_{q_{i}} $.
For example, for any practical purpose, the query $ q = $``R.A. Fisher was a statistician. Was he great?" is an element of $ \mathcal{Q}_{0} $ for $ q_{0} = $``RA Fisher was a statistician. Was he great?". 
Defining $ \mathcal{F}_{0} := \{F_{q}: q \in \mathcal{Q}_{0} \} $, 
we modify  Eq.\ \eqref{eq:primary-interest} to  
\begin{equation}
\label{eq:new-approach}
    H_{0}:  \min_{F_{q}\in \mathcal{F}_{0}} d(F_{q'}, F_{q})=0 \quad \quad \quad vs \quad \quad \quad H_{A}:  \min_{F_q \in \mathcal{F}_{0}} d(F_{q'}, F_{q})>0 
\end{equation} for some distance $ d $ defined for distributions on $ \mathbb{R}^{s} $.

As with the test described in Eq. \eqref{eq:primary-interest}, we assume that the map from the space of queries to their corresponding response distributions is not known and we must obtain an estimate of each $ F_{q} $, denoted by $\hat{F}_{q}$, given responses $ f(q)_{1}, \hdots, f(q)_{r} $ for each $ q \in \mathcal{Q}_{0} \cup \{ q'\}$.
In practice, repeatedly querying $ f $ may be prohibitively expensive, especially for large $ | \mathcal{Q}_{0} | $. 
The remainder of this paper describes and analyzes a hypothesis test of the form described in Eq. \eqref{eq:new-approach} in the setting where the response space is restricted to $ \{0, 1\} $ and the user has a known resource budget $ \nu $.

\paragraph{Notation.} For any $x \in \mathbb{R}$, $\lceil x \rceil$
denotes ceiling -- the smallest integer just exceeding $x$, and $\lfloor x \rfloor$ denotes floor -- the largest integer not exceeding $x$. For any natural number $n \in \mathbb{N}$, $[n] \coloneq \lbrace 1,2,\dots ,n \rbrace$. For a set of values $p_1,\dots,p_m$, the order statistics are defined as $p_{(1)} \leq p_{(2)} \leq \dots \leq p_{(m)}$, that is, $p_{(i)}$ is the $i$-th minimum value in the set $\lbrace p_1,\dots,p_m \rbrace$.




\section{Preliminaries}
\label{sec:preliminaries}

\subsection{Basics of the statistical hypothesis testing framework}

Suppose we observe a sample $ X_1, \dots, X_r \sim^{i.i.d.} F $  where the distribution function $ F $ is parameterized by $ \theta \equiv \theta(F) \in \Theta $.
The goal of parametric statistical hypothesis testing is to determine if $ \theta $ is an element of the set $ \Theta_{0} \subset \Theta $ or an element of the set $ \Theta_{1} = \Theta \setminus \Theta_0 $ 
based on a test statistic $ T_{r} \equiv T_{r}(X_{1}, \hdots, X_{r}) $. 
In the statistical hypothesis testing framework, we reject $H_0: \theta \in \Theta_{0} $  
in favor of the alternative hypothesis $H_A:\theta \in \Theta_1$
only if the observed sample provides sufficient evidence against it. Denoting the set of all possible values of the test statistic $T_r$ by $\mathcal{T}^r$, we select a rejection region $\mathcal{T}^r_1 \subset \mathcal{T}^r$ such that we reject $H_0$ if and only if the observed value of the test statistic $t_r \equiv t_{r}(x_{1}, \hdots, x_{r}) \in \mathcal{T}^r_1$.

There are two types of errors in this setting: Type-I and Type-II. 
A Type-I error is the rejection of $H_0$ when it is true; that is, when $\theta \in \Theta_0$ but the practitioner determines $\theta \in \Theta_1$ based on $t_r \in \mathcal{T}^r_A$.
A Type-II error is the failure to reject $H_0$ when it is false; that is, when $\theta \in \Theta_1$ but the practitioner determines $\theta \in \Theta_0$ based on $t_r \notin \mathcal{T}^r_1$.
For a given test statistic $ T_{r} $, there is an inherent tradeoff between $\mathbb{P}[\text{Type-I error}]$ and $\mathbb{P}[\text{Type-II error}]$. 
In the Neyman-Pearson framework, the user defines a tolerance level for $\mathbb{P}[\text{Type-I error}]$, and then selects a $ T_{r} $ that minimizes $\mathbb{P}[\text{Type-II error}]$. 

Formally, a testing procedure is a binary function $\gamma:\mathcal{T}^{r} \to \{0,1\}$ such that 
\begin{equation*}
\begin{aligned}
\gamma_r \equiv
\gamma(t_r)=
\begin{cases}
    1 \;\; \text{ if } t_r \in \mathcal{T}^{r}_{1}   \\
      0 \;\; \text{ otherwise}.
      \\
\end{cases}
\end{aligned}
\end{equation*} 

A few important definitions are given below.
\begin{definition}
    The size of a testing procedure $ \gamma_r $ is defined as 
    \begin{equation*}
        \mathrm{size}(\gamma_r)=\sup_{\theta \in \Theta_0}
        \mathbb{P}_{\theta}[\gamma_r(T_{r}) = 1]. 
    \end{equation*}
\end{definition}



As mentioned above, in the Neyman-Pearson framework the user specifies a Type-I error tolerance (or level of significance, denoted $\alpha$) and considers only testing procedures for which the probability of making a Type-I error is controlled as specified. 
Tests with this property are referred to as {\em valid}.
\begin{definition}
    Given a level of significance $\alpha$, the testing procedure $ \gamma $ is valid if it has size less than or equal to $\alpha$; that is, 
    $\sup_{\theta \in \Theta_0}\mathbb{P}[\gamma_r(T_{r}) = 1] \leq \alpha$. 
\end{definition}

In practice, some testing procedures 
may not be valid for a given $ r $ but approach validity as $r$ grows. 
These procedures might have other desirable properties which warrant their use despite them not being strictly valid, and are termed \textit{asymptotically valid} tests. 
\begin{definition}
    Given a level of significance $\alpha$, a sequence of testing procedures $ \lbrace \gamma_r \rbrace_{r=1}^{\infty}$ is asymptotically valid if their sizes approach $\alpha$, that is, 
    \begin{equation*}
    \lim_{r \to \infty}
        \mathrm{size}(\gamma_r)
  \leq \alpha.  \end{equation*}
\end{definition}

Given an 
(asymptotically) valid testing procedure, the remaining consideration is the probability of committing a Type-II error. 
For this,  the {\em power function} is introduced. 
\begin{definition}
    The power function $\beta_{\gamma_r}:\Theta_1 \to [0,1]$ of a testing procedure $ \gamma_r$ is the probability of rejecting $ H_0 $ as a function of the parameter value $ \theta $; that is,
    \begin{equation*}
        \beta_{\gamma_r}(\theta)=\mathbb{P}_{\theta}[\gamma_r(T_{r}) = 1].
    \end{equation*}
\end{definition}

A desirable property of an (asympotically) valid testing procedure is the power function approaching $ 1 $ for all $ \theta \in \Theta_{1} $ as $ r $ grows.
A testing procedure with this property is called {\em consistent}.

\begin{definition}
   A sequence of testing procedures $ (\gamma_{1}, \hdots, \gamma_{r}) $ is consistent if the sequence of power functions $ (\beta_{\gamma_{1}}, \hdots, \beta_{\gamma_{r}} ) $ approaches $ 1 $ for all $ \theta \in \Theta_{1} $; that is
   \begin{equation*}
       \lim_{r \to \infty}
       \beta_{\gamma_{r}}(\theta)=1 \;\;\;\;\text{for all }\theta \in \Theta_1.
   \end{equation*}
\end{definition}

A testing procedure that is (asymptotically) valid and consistent properly controls the Type-I error and correctly rejects $ H_{0} $ as $ r $ grows, the two most fundamental properties within the Neyman-Pearson framework. The key technical contribution of this paper is to devise an asymptotically valid test which controls the number of the undesirable rejections demonstrated in Figure \ref{fig:motivation}, under realistic budget constraints. 

\subsection{Statistical methods for generative models}
The recent improvements in the accessibility and performance of generative models for everyday uses \citep{jiang2023clip, grattafiori2024llama, achiam2023gpt, anthropic2024claude3family, ustun2024aya,team2023gemini}
and similar improvements in specialized domains such as medicine \citep{thirunavukarasu2023large, nori2023can, abd2023large}, radiology \citep{d2024large, kim2024large}, law \citep{sun2023short, 10850911}, etc. \citep{lo2023impact,rahman2023chatgpt,helm2023statistical,khan2024chatgpt, zhang2024codeagent}, has spurred investigations into failure modes of the models and the systems in which they are embedded.
For example, \cite{ness2024medfuzz} demonstrated that model performance on a popular medical benchmark is highly sensitive to medically-irrelevant insertions and perturbations; \cite{gallifant2024language} showed models may be over reliant on knowledge of the name of a drug as opposed to its properties;  \cite{chen2024premise} showed that the order of independent premises in a logical statement can affect performance by up to 30\%.

The demonstration of simple but pervasive failure modes has motivated the application and development of statistical methods for understanding and comparing relevant properties of models, conditionings, and prompt structures.
Different applications and different methods require different model accessibility assumptions -- for example, it is possible to determine if a model was trained on a particular type of data when token-wise log-probabilities of a relevant set of tokens are available \citep{shi2023detecting} as well as when users only have access to model responses \citep{helm-etal-2025-statistical}.
Given that access to the token-wise log-probabilities or other model internals also implies access to the model responses, we choose to operate in the setting of access only to the responses.
Perhaps the most fundamental statistical development in this paradigm is the treatment of model evaluation as a statistical problem that requires comparing distributional properties of outputs and scores before making declarative statements \citep{miller2024adding}.
The current paper builds on this treatment of model evaluation by developing a statistical hypothesis test that addresses the wide class of aforementioned observed failure modes and is asymptotically valid.

\subsection{Testing of unspecified null hypothesis}
In the traditional hypothesis testing framework, when testing  $H_0: \theta' \in \Theta_0$, the null region $\Theta_0$ is specified. However, in our case, this problem is extended to testing $H_0:\theta' \in \Theta_0$ where $\Theta_0$ is unspecified. In such case, 
one can assume that the practitioner has the ability to draw a random sample $\theta_1,\dots,\theta_m \in \Theta_0$. If this null sample has sufficient coverage of $\Theta_0$, a large deviation of $\theta'$ from 
$\lbrace \theta_1,\dots,\theta_m \rbrace$ provides evidence against $H_0$. 

In our case we do not observe the $\theta_i$ directly, but we can sample from their distributions; thus, in addition to drawing samples from distribution $F_{\theta'}$ (a probability distribution characterized by $\theta$) to compute $\hat{\theta'}$,  random samples are also drawn from $F_{\theta_1},\dots,F_{\theta_m}$ to compute estimates $\hat{\theta}_1,\dots, \hat{\theta}_m$. A sufficiently large deviation of the estimate $\hat{\theta'}$ from the estimated null sample $\lbrace \hat{\theta}_1,\dots,\hat{\theta}_m \rbrace$, indicating deviation of $\theta'$ from the null sample $\lbrace \theta_1,\dots,\theta_m \rbrace$, leads to the rejection of the null hypothesis $H_0:\theta' \in \Theta_0$, if the null sample has sufficient coverage of $\Theta_0$.
A similar technique has been used to empirically calibrate p-values in observational studies for drug-safety \citep{schuemie2014interpreting}. In our paper, the set of Bernoulli parameters of all possible semantically irrelevant perturbations to $q_0$, denoted $\mathcal{P}_0$, is unknown.

\subsection{Stability of statistical results to reasonable perturbations}
Our investigation is motivated by the observation that a conventional statistical hypothesis test often concludes that significant change in response distribution has occurred due to a semantically irrelevant perturbation to a query. In \cite{yu2024veridical} and \cite{agarwal2025pcs}, the authors discuss the principle of stability of statistical results relative to reasonable perturbations in the data and the model, which makes statistical results reproducible. 
Our work on the case of a composite null hypothesis consisting of semantically irrelevant perturbations to a query is an example of investigating stability of the test to ``reasonable perturbations'' to the null query. 
\section{Methodology}
\label{Sec:Model_and_Method}

As in our motivating example, we focus our analysis in the setting where $ \mathcal{Q} $ is restricted to queries where $ F_{q} $ has two elements in its support and, thus, $ F_{q} $ is completely parameterized by a Bernoulli parameter $ p $.
Our goal is to test if the Bernoulli parameter of a test query $q'$, denoted by $p'$, is close to the Bernoulli parameter of a null query $q_0$, denoted by $p_0$, while taking into account a user-defined notion of semantic similarity.
In particular, let $ \mathcal{Q}_{0} $ be a set of queries deemed semantically equivalent to $ q_{0} $, and $ \mathcal{P}_{0} $ be the set of corresponding Bernoulli parameters.
Thus, our goal is to test $H_0:p' \in \mathcal{P}_0$ against $H_A:p' \notin \mathcal{P}_0$.

Since the mapping from $ \mathcal{Q}_{0} $ to $ \mathcal{P}_{0} $ is not known \textit{a priori}, we must estimate some of the elements of $ \mathcal{P}_{0} $.
We first sample queries $q_1,\dots,q_m \in \mathcal{Q}_0$. For every sampled query $q_j$, we obtain $i.i.d.$ responses $f(q_j)_1,\dots,f(q_j)_r \sim^{i.i.d.} \mathrm{Bernoulli}(p_j)$, and
estimate the Bernoulli parameter $p_j$ by $\hat{p}_j=\frac{1}{r} \sum_{k=1}^r f(q_j)_k$. 
Similarly, for the test query $q'$, we obtain $i.i.d.$ responses $f(q')_1,\dots,f(q')_r \sim^{i.i.d.} \mathrm{Bernoulli}(p')$, and estimate the 
Bernoulli parameter $p'$ by $\hat{p}'=\frac{1}{r} \sum_{k=1}^r f(q')_k$.
Subsequently, we define our test statistic as
\begin{equation}
\label{Eq:realistic_test}
    T_{m,r}=\min_{j \in [m]}|\hat{p}_j-\hat{p}'|
\end{equation}
and reject $H_0:p' \in \mathcal{P}_0$ if $T_{m,r}>\epsilon$ for an appropriately chosen $\epsilon$. 

For given test query $q'$, and choices for $m$ and $r$, we provide the procedure for computing $T_{m,r}$ in Algorithm \ref{Algo:TestGeneric}. 
\begin{algorithm}[h!]
\caption{GenericStatistic($f$, $ \mathcal{Q}_0 $, $q'$, $m$, $ r$)} 
\label{Algo:TestGeneric}
\begin{algorithmic}[1]
\STATE Sample
$i.i.d.$ queries $q_{1}, q_{2},\dots q_{m} \in \mathcal{Q}_0 $.
\FOR{$j \in \{1, \hdots, m\}$}
    \STATE Sample $i.i.d.$ replicates $ f(q_{j})_{1}, \hdots, f(q_{j})_{r} $.
    \STATE $ \hat{p}_{j} \leftarrow \frac{1}{r}\sum_{k=1}^{r} f(q_{j})_{k} $.
\ENDFOR \\
\STATE Sample $i.i.d.$ replicates $ f(q')_{1}, \hdots, f(q')_{r} $.
\STATE $ \hat{p}' \leftarrow \frac{1}{r}\sum_{k=1}^{r} f(q')_{k} $.
\STATE $T_{m,r} \leftarrow \min_{j \in \{1, \hdots, m\}}
\left| \hat{p}_j-\hat{p}' \right|$.
\RETURN $T_{m,r}$.
\end{algorithmic}
\end{algorithm}

We note that the size and power function of the test based on $T_{m,r}$,
depend on $\epsilon, m, $ and $ r $. 
Moreover, for any fixed $ r $, increasing  $ m $ decreases the power of the test, because, even when $H_A:p' \notin \mathcal{P}_0$ is true, the probability of at least one estimate $\hat{p}_j$ behaving erratically (i.e., is far from $ p_{j} $) and being close to $\hat{p}'$ increases, thereby decreasing the probability of rejecting $H_0$. 
As such, there is a natural interplay between $ m, r,$  and $ \epsilon $ that affects the properties of the proposed test. The question, then, is how to choose $ m $, $ r $, and $ \epsilon $, given level of significance $\alpha$ and budget $\nu$.  

\subsection{Proposed test under realistic budget constraint}
\label{Subsec:Optimal_Test}
We extend $\mathcal{P}_0$ to the interval $[a,b]=[\min_{p \in \mathcal{P}_0}, \max_{p \in \mathcal{P}_0}]$. 
In reality, we operate under a budget constraint $m \cdot r \leq \nu$ and do not know the parameters $a$ and $b$.
We compute estimates $\hat{a}$ and $\hat{b}$ with $ \tilde{m}$, $ \tilde{r} $ such that $ \tilde{m} \cdot \tilde{r} \ll \nu $, via the procedure is described in Algorithm \ref{Algo:EstimateRange}. 

\begin{algorithm}[h!]
\caption{EstimateRange($f$, $\mathcal{Q}_0$, $\tilde{m}$, $\tilde{r}$)}. 
\label{Algo:EstimateRange}
\begin{algorithmic}[1]
\STATE Generate null queries $q_1,\dots,q_{\tilde{m}} \in \mathcal{Q}_0$. 
\FOR{$ j \in \lbrace 1, 2,  \hdots, \tilde{m}  \rbrace$}
\STATE Obtain $i.i.d.$ responses $f(q_j)_1,\dots,f(q_j)_r$. 
\STATE $ \hat{p}_j \leftarrow \frac{1}{\tilde{r}} \sum_{k=1}^{\tilde{r}} f(q_j)_k$. 
\ENDFOR
\STATE $\hat{a} \leftarrow \hat{p}_{(1)}= \min_{j \in [\tilde{m}]} \hat{p}_j$, $\hat{b} \leftarrow \hat{p}_{(\tilde{m})}=\max_{j \in [\tilde{m}]} \hat{p}_j$. 
\RETURN $(\hat{a},\hat{b})$.
\end{algorithmic}
\end{algorithm}

 In Section \ref{Sec:Theoretical_results_h}, under a set of technical assumptions,
we derive the expressions for a validity constraint (in Theorem \ref{Thm:lower_bound_on_real_size_h}) and a lower bound on average power (in Theorem \ref{Thm:average_power_real_test_h}), which involve the unknown parameters $a$ and $b$.
Hence, we approximate the expressions using the output $(\hat{a},\hat{b})$ of Algorithm \ref{Algo:EstimateRange}. We choose the triple $(\epsilon^{**},m^{**},r^{**})$ by maximizing the said approximate lower bound on average power, given by
\begin{equation*}
    \hat{H}(\epsilon,m,r)=
    \frac{2}{\hat{b}-\hat{a}}
\left\lbrace
    \left(
1-
\frac
{
\epsilon+
\sqrt{\frac{\mathrm{log}r}{r}}
}
{
\hat{b}-\hat{a}
}
    \right)^m -1 
\right\rbrace
\left(
\epsilon+ 
\sqrt{\frac{\mathrm{log}r}{r}}
\right)+ 
\left(
1-\frac{2m}{\sqrt{r}}
\right),
\end{equation*}
where $(\epsilon,m,r)$ are such that they satisfy the approximate validity constraint
\begin{equation*}
    1- 
    \frac{1}{\hat{b}-\hat{a}}
    \left(
\epsilon-
\sqrt{\frac{\mathrm{log}r}{r}}
    \right)^m + 
    \frac{2m}{\sqrt{r}}
    \leq \alpha.
\end{equation*}
The entire testing procedure is described in Algorithm \ref{alg:test-optimal}.
\newline
\begin{algorithm}[h!]
\caption{OptimalTest($f$, $\mathcal{Q}_0$, $ q' $,
$\alpha$,
$\nu$,
 $ \tilde{m} $, $ \tilde{r}$,  $\eta_{\epsilon} $, $ \epsilon_{\max} $)}. 
\label{Algo:TestOptimal}
\begin{algorithmic}[1]
\STATE $ \hat{a}, \hat{b} \leftarrow$ EstimateRange($f$,$\mathcal{Q}_0,\tilde{m}$, $\tilde{r}$).
\STATE $ \nu \leftarrow \nu - \tilde{m} \cdot \tilde{r} $.
\STATE $\epsilon_{\max} \leftarrow \min
\left\lbrace 
\hat{a},\hat{b}-\hat{a},1-\hat{b}
\right\rbrace
$.
\STATE $ h^{*} \leftarrow -\infty $,
$\epsilon^{**} \leftarrow 0$, $m^{**} \leftarrow 1$, $r^{**} \leftarrow 1$. 
\FOR{$ \epsilon \in [0, \eta_{\epsilon} $, $2 \eta_{\epsilon}, \hdots,  \epsilon_{\max}]$}
\STATE $ m \leftarrow
\max
\left\lbrace
\lceil 
\frac{
|\mathrm{log}(\alpha)|
}
{
|\mathrm{log}(1-\frac{\epsilon}{\hat{b} - \hat{a}})|
}
    \rceil, \tilde{m}
\right\rbrace 
    $.
\STATE $ r \leftarrow \frac{\nu}{m} $.
\STATE \texttt{is\_valid} $ \leftarrow \mathbbm{1}\left\{\left(
1-
\frac{1}{\hat{b}-\hat{a}}
\left(
\epsilon-
\sqrt{\frac{\mathrm{log}r}{r}}
\right)
\right)^{m}
+\frac{2m}{r}
\leq \alpha\right\} $. 
\IF{\texttt{is\_valid}}
\STATE $ \hat{h} \leftarrow \frac{2}{1-(\hat{b}-\hat{a})}
    \left\lbrace
\left(
1-
\frac
{
\epsilon+
\sqrt{\frac{\mathrm{log}r}{r}}
}{\hat{b}-\hat{a}} 
\right)^{m} -
1
    \right\rbrace
    \left(
\epsilon+
\sqrt{\frac{\mathrm{log}r}{r}}
    \right) +
    \left(
1- \frac{2m}{\sqrt{r}}
    \right) $
\IF{$ \hat{h} > h^{*} $}
\STATE $ h^{*}\leftarrow \hat{h}$,
$\epsilon^{**} \leftarrow \epsilon$,
$m^{**} \leftarrow m$,
$r^{**} \leftarrow r$.
\ENDIF
\ENDIF
\ENDFOR
\STATE $ T \leftarrow  \mathrm{GenericStatistic}$($f$,$\mathcal{Q}_0 $, $q'$, $m^{**}$, $ r^{**}$).
\RETURN $ \mathbbm{1}\{T > \epsilon^{**}\} $
\end{algorithmic}
\label{alg:test-optimal}
\end{algorithm}
\begin{remark}
    Under the condition that the Bernoulli parameters $p_1,\dots,p_m \sim^{i.i.d.} \mathrm{Unif}([a,b])$, one can set
    $$\hat{a}=
    \left(
    \hat{p}_{(1)}-
    \frac{1}{\tilde{m}+1}
    (\hat{p}_{(\tilde{m})}-\hat{p}_{(1)})
    \right)
    , 
    \hat{b}=
    \frac{\tilde{m}+1}{\tilde{m}}
    \hat{p}_{(\tilde{m})}
    $$
    in Algorithm \ref{Algo:EstimateRange}, for bias correction. 
\end{remark}
\section{Theoretical Results}
\label{Sec:Theoretical_results_h}
We state our theoretical results in this section. 
We first briefly recall our setting once again.
\newline
We have a generative model $f$ which provides binary responses (``Yes'',``No'') to any query, so that response to any query $q$ can be treated as a Bernoulli random variable. Let $\mathcal{Q}$ be the set of all queries and let $q_0 \in \mathcal{Q}$ be a particular query. We define $\mathcal{Q}_0$ to be the set of queries which are semantically similar to $q_0$. Let $\mathcal{P}_0=[a,b]$ denote the smallest interval containing the set of all Bernoulli parameters corresponding to the queries in $\mathcal{Q}_0$;  $a$ and $b$ are unknown. For a new query $q' \in \mathcal{Q}$, suppose the corresponding Bernoulli parameter is  $p'$. We want to test $H_0: p' \in \mathcal{P}_0$ against $H_A:p' \notin \mathcal{P}_0$. However, since the parameters $a$ and $b$ are unknown, we adopt the following strategy. We generate queries $q_1,\dots,q_m \in \mathcal{Q}_0$, estimate their Bernoulli parameters, and reject $H_0$ if the estimated Bernoulli parameter of the test query $q'$ is sufficiently far from  the estimated Bernoulli  parameters of each of the sampled null queries $q_1,\dots,q_m$. For any query $q$, we estimate the corresponding Bernoulli parameter by the mean of $r$ Bernoulli responses to that query.

\begin{figure}[t!]
    \centering
    \includegraphics[width=\linewidth]{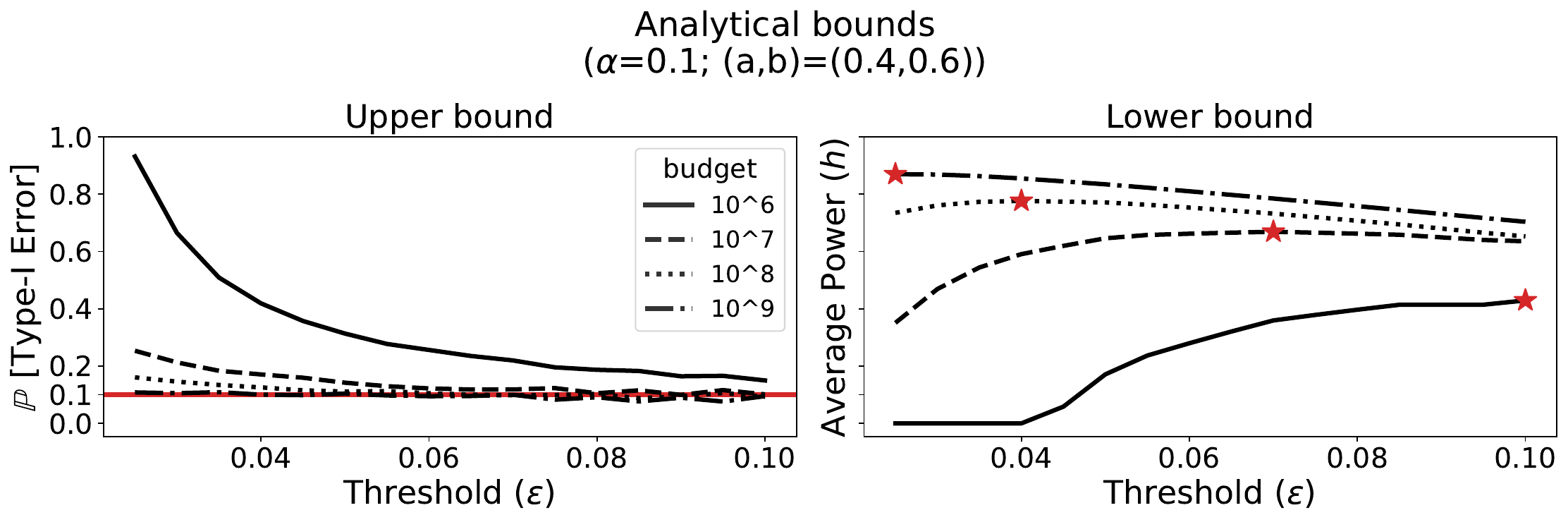}
    \caption{Analytical upper bounds for the Type I Error (left) and lower bounds for the average power (right) of the proposed test for various thresholds ($ \epsilon$) and budgets.
    The maximal average power for a given budget is highlighted by a red star.
     Algorithm \ref{alg:test-optimal} selects $ m, r, $ and $ \epsilon $ such that the test is asymptotically valid ($\mathbb{P}[\text{type-I error}] \to \alpha$), and an estimated lower bound on the average power is maximized.}
    \label{fig:analytical-bounds}
\end{figure}

\subsection{The ideal test}

To develop a method for choosing $ m, r, $ and $ \epsilon $ we introduce an ideal version of Eq. (\ref{Eq:realistic_test}) where the Bernoulli parameters  $p_1,\dots,p_m,p'$ are known. 
We define
\begin{equation}
\label{Eq:ideal_test}
\tilde{T}_{m} = \min_{j \in [m]} |p' - p_j |,
\end{equation} 
and reject $ H_{0} $ if $ \tilde{T}_{m} > \epsilon $ for a suitably chosen $ \epsilon > 0 $.
We refer to this test procedure based on $ \tilde{T}_{m} $ as the ``ideal" test and the test procedure based on $ T_{m,r} $ as the ``realistic" test.
We derive analytical upper bound for the size and lower bound for the power of the ideal test, and show that the realistic and ideal tests are close for a sufficiently large budget, and choose $ m, r, $ and $ \epsilon $ based on corresponding bounds on the realistic test.
We finalize our theoretical setting with the following three assumptions.
\begin{assumption}
\label{Asm:null_alt_independence_h}
The random vector $(p_1,\dots,p_m)$ is independent of $p'$, that is, their joint pdf can be written as
\begin{equation*}
    f(p_1,\dots,p_m,p')=f(p_1,\dots,p_m)f(p').
\end{equation*}
\end{assumption}
Assumption \ref{Asm:null_alt_independence_h} ensures that test statistics for the ideal test and the realistic tests are close. 
\begin{assumption}
\label{Asm:null_uni_dist_h}
    For queries $ q_{1}, \hdots, q_{m} $ randomly sampled from $ \mathcal{Q}_{0} $, the corresponding Bernoulli parameters are uniformly distributed on $ [a,b]$; that is,
    \begin{equation*}
        p_1,\dots,p_{m} \sim^{i.i.d.} \mathrm{Unif}[a,b] 
    \end{equation*}
    where $0<a<b<1$. 
\end{assumption}
Assumption \ref{Asm:null_uni_dist_h} is needed to establish an upper bound on the sizes of the ideal test and the realistic tests. 
\begin{assumption}
\label{Asm:alt_uni_dist_h}
    For a query $ q' $ randomly sampled from $ \mathcal{Q} $, the corresponding Bernoulli parameter is uniformly distributed on $ (0,1) $; that is, 
    \begin{equation*}
        p' \sim \mathrm{Unif}(0,1).
    \end{equation*}
\end{assumption}
Assumption \ref{Asm:alt_uni_dist_h} allows us to obtain an approximate lower bound on the power function of the realistic test.

Our first result shows that for any fixed $(p_1,\dots,p_m,p')$, the difference between $ T_{m,r} $ and $ \tilde{T}_{m} $ approaches zero as $m,r \to \infty$, if $r=\omega(m^2)$. Thus, if the budget $ \nu = m \cdot r \to \infty $ such that $ r=\omega(m^2) $, we can approximate $\tilde{T}_m$ with $T_{m,r}$.

\noindent 

\begin{lemma}
\label{lemma:concentration_of_test_statistic}
Suppose that for every query $ q \in \{q_{1}, \hdots, q_{m}, q'\} $ we observe $i.i.d.$ replicates of responses denoted by
$f(q)_1,\dots,
f(q)_r \sim^{iid} \mathrm{Bernoulli}(p)$  where $p$ is the Bernoulli parameter of the query $q$. 
Define 
\begin{equation*}
\begin{aligned}
    &\hat{p}_j=\frac{1}{r} \sum_{k=1}^r f(q_j)_k \text{ for all $j \in [m]$}, \\
    &\hat{p}'=\frac{1}{r} \sum_{k=1}^r f(q')_k, 
\end{aligned}
\end{equation*}
and set $\tilde{T}_m=\min_{j \in [m]}|p_j-p'|$ and $T_{m,r}=\min_{j \in [m]} |\hat{p}_j-\hat{p}'|$. 
Then, conditioning on $ (p_{1},\dots p_{m}, p')$,
\begin{equation*}
\mathbb{P}
\left[
   \big|
T_{m,r}-\tilde{T}_m
   \big|
   <
   \sqrt{\frac{\mathrm{log}r}{r}}
   \bigg|
   (p_1,\dots p_{m},p')
\right] \geq 1 -
\frac{2m}{\sqrt{r}}
.
\end{equation*}
\end{lemma}
 
The proof of Lemma \ref{lemma:concentration_of_test_statistic} is provided in the Appendix.
An important extension of Lemma \ref{lemma:concentration_of_test_statistic} is a bound for the difference conditioning only on $ p' $.

\begin{theorem}
\label{thm:concentration_of_test_statistic_only_ew_conditioning}
Suppose Assumption \ref{Asm:null_alt_independence_h} holds, and consider 
the setting of \textit{Lemma \ref{lemma:concentration_of_test_statistic}}. For any $p' \in (0,1)$,
\begin{equation*}
    \mathbb{P}
    \left[
    \left|
T_{m,r}-\tilde{T}_m
    \right|
    < 
    \sqrt{\frac{\mathrm{log}r}{r}}
    \bigg| p'
    \right]
    \geq 
1- \frac{2m}{\sqrt{r}}
    . 
\end{equation*}
\end{theorem}

The proof of Theorem \ref{thm:concentration_of_test_statistic_only_ew_conditioning} is provided in the Appendix and is based on the independence of $ p' $ and $ p_{1}, \hdots, p_{m} $.
Theorem \ref{thm:concentration_of_test_statistic_only_ew_conditioning} ensures the two test statistics are close for all possible samples from $ \mathcal{Q}_{0} $.

\subsection{Control over size}
\label{Subsec:Control_size}
We next provide bounds for the size and power of the ideal and realistic tests.
We start by deriving a lower bound on the number of queries required for the ideal test to be valid.

\begin{lemma}
\label{lemma:m_for_true_valid}
Under \textit{Assumption \ref{Asm:null_alt_independence_h}} and \textit{Assumption \ref{Asm:null_uni_dist_h}}  
on $ q_{1}, \hdots, q_{m} $, suppose we observe the true corresponding Bernoulli parameters $ p_{1}, \hdots, p_{m}, p' $ and we reject $H_0$ if $\tilde{T}_m>\epsilon$.
Then, if $\epsilon<\min \lbrace a,b-a,1-b \rbrace$, a sufficient condition to ensure that the test is valid at level of significance $\alpha$ is given by
\begin{equation}
\label{eq:m}
    m \geq 
    \Bigl\lceil
    \frac
    {
|\mathrm{log}(\alpha)|
    }
    {
    |\mathrm{log}(1-\frac{\epsilon}{b-a})|
    }
    \Bigl\rceil
    .
\end{equation}
\end{lemma}

The closeness of the ideal test to the realistic test suggests that a sufficient condition for validity of the realistic test is close to this sufficient condition for validity of the ideal test. 
However, since we have to estimate the Bernoulli parameters in the realistic test, apart from ensuring $m$ is sufficiently large, we also need to ensure $r$ is sufficiently large. 
This is established by our next result, which says that if $m$ and $r$ are large enough, then the realistic test is valid. 
\begin{theorem}
\label{Thm:lower_bound_on_real_size_h}
Suppose 
\textit{Assumption \ref{Asm:null_alt_independence_h}} and
\textit{Assumption \ref{Asm:null_uni_dist_h}} hold, and the Bernoulli parameters $p_1,\dots,p_m$ and $p'$ are not observed. For every query $q \in \lbrace
q_1,\dots,q_m,q'
\rbrace$, $i.i.d.$ replicates of responses, denoted by $f(q)_1,\dots,f(q)_r$ are obtained. Define $\hat{p}_j=\frac{1}{r} \sum_{k=1}^r f(q_j)_k$ and $\hat{p}'=\frac{1}{r} \sum_{k=1}^r f(q')_k$. 
Recall that
for testing $H_0:p' \in \mathcal{P}_0$ versus $H_A:p' \notin \mathcal{P}_0$,
our decision rule rejects $H_0$ if $T_{m,r}>\epsilon$, for a chosen threshold $\epsilon$, where 
$T_{m,r}=\min_{j \in [m]}|\hat{p}_j-\hat{p}'|$. 
 When $\epsilon<\min \lbrace a,b-a,1-b \rbrace$ and
$r$ is sufficiently large, for all $p' \in \mathcal{P}_0$, 
\begin{equation*}
\begin{aligned}
    \mathbb{P}
    \left[
    T_{m,r}>\epsilon
    \bigg| p'
    \right] 
    \leq 
    \left(
1- 
\frac
{
\epsilon -
\sqrt{\frac{\mathrm{log}r}{r}}
}{b-a}
    \right)^m +
    \frac{2m}{\sqrt{r}}.
\end{aligned}
\end{equation*}
\end{theorem}



Note that the upper bound  on the size of the realistic test approaches zero as $m \to \infty$ and $r \to \infty$ such that $r=\omega(m^2)$.
We plot the behavior of the upper bound on the size for the realistic test for $ (a,b) = (0.4, 0.6) $ and $ \alpha = 0.1 $ as a function of $ \epsilon $ for various budgets in the left panel of Figure \ref{fig:analytical-bounds}.

\subsection{Control over power}
\label{Subsec:Control_power}
In this subsection, we establish lower bounds for power of the ideal and realistic tests. 
First, we state a lemma which provides an expression for the power function of the ideal test.
\begin{lemma}
\label{Lm:power_for_true_test}
In our setting, under \textit{Assumptions} 
\ref{Asm:null_alt_independence_h} and
\ref{Asm:null_uni_dist_h}, suppose we observe the true Bernoulli parameters, and we want to test $H_0:p' \in \mathcal{P}_0$ versus $p' \notin \mathcal{P}_0$ at level of significance $\alpha$. Our ideal decision rule rejects $H_0$ if $\tilde{T}_m>\epsilon$ for some chosen threshold $\epsilon$, where $\tilde{T}_m=\min_{j \in [m]} |p_j-p'|$. 
For $\epsilon<\min \lbrace a,b-a,1-b \rbrace$, the power function of the ideal test is given by
\begin{equation*}
\begin{aligned}
\tilde{\beta}(p')
\equiv \tilde{\beta}_{m,\epsilon}(p')
\coloneq
    \mathbb{P}
    \left[
    \tilde{T}_m>\epsilon
    \bigg| p'
    \right]=
    \left\{
\begin{array}{ll}
      1,
       & p'\in (0,a-\epsilon] \cup [b+\epsilon,1) \\
       \left(
      \frac{b-p'-\epsilon}{b-a}
      \right)^m,
       & p' \in (a-\epsilon,a) \\
       \left(
      \frac{p'-\epsilon-a}{b-a}
      \right)^m, & p' \in (b,b+\epsilon)  \\
\end{array} 
\right.
\end{aligned}
\end{equation*}
\end{lemma}
With the help of Theorem \ref{thm:concentration_of_test_statistic_only_ew_conditioning} and the abovementioned Lemma \ref{Lm:power_for_true_test}, we deduce a lower bound on the power of the realistic test. 
\begin{theorem}
\label{Thm:lower_bound_power_realistic_h}
Consider the setting of \textit{Theorem \ref{Thm:lower_bound_on_real_size_h}}. 
Recall that for testing $H_0:p' \in \mathcal{P}_0$ versus $H_1:p' \notin \mathcal{P}_0$,
our realistic decision rule rejects $H_0$ if $T_{m,r}>\epsilon$, for a chosen threshold $\epsilon$, where 
$T_{m,r}=\min_{j \in [m]}|\hat{p}_j-\hat{p}'|$. 
When $\epsilon<\min \lbrace a,b-a,1-b \rbrace$ and $r$ is sufficiently large,
$\mathbb{P}[T_{m,r}>\epsilon|p'] \geq \phi(p')$ where the lower bound $\phi$ is given by,
\begin{equation*}
\begin{aligned}
    \phi(p')
    \coloneq
    \left\{
\begin{array}{ll}
    \left(
1-\frac{2m}{\sqrt{r}}
    \right)
      ,
       & p'\in \left(0,a-\epsilon-
    \sqrt{\frac{\mathrm{log}r}{r}}
       \right]  \cup \left[
       b+\epsilon+
    \sqrt{\frac{\mathrm{log}r}{r}}
       ,1 
       \right) \\
       \left(1-
      \frac{\epsilon+
    \sqrt{\frac{\mathrm{log}r}{r}}
      }{b-a}
      \right)^m
      - \frac{2m}{\sqrt{r}}
      ,
       & p' \in (a-\epsilon-
    \sqrt{\frac{\mathrm{log}r}{r}}
       ,a) \\
       \left(1-
      \frac{\epsilon+
    \sqrt{\frac{\mathrm{log}r}{r}}
      }{b-a}
      \right)^m
      - \frac{2m}{\sqrt{r}}
      , & p' \in (b,b+\epsilon+
    \sqrt{\frac{\mathrm{log}r}{r}}
      )  \\
\end{array} 
\right.
\end{aligned}
\end{equation*}
for all $p'\in \mathcal{P}_1 \coloneq (0,1) \setminus \mathcal{P}_0= (0,a) \cup (b,1)$,
\end{theorem}

Based on the results established in this section, we derive a sufficient condition for consistency and asymptotic validity, which is stated below. 



\begin{corollary}
\label{Cor:asy_valid_const_true_par_h}
Consider the setting of \textit{Theorem \ref{Thm:lower_bound_power_realistic_h}}. As $\epsilon \to 0$, $m,r \to \infty$ such that $r=\omega(m^2)$,
$\sqrt{\frac{\mathrm{log}r}{r}} \leq \epsilon$ and $\epsilon-\sqrt{\frac{\mathrm{log}r}{r}} \leq (b-a)$,
 we have an  asymptotically valid and consistent sequence of tests. 
\end{corollary}

\subsection{Choosing $\epsilon$, $ m $ and $ r $}
We first ensure our chosen $\epsilon,m$ and $r$ approximately satisfy the validity constraint. Amongst the selected values of $\epsilon,m,r$, we intend to choose those which maximize the average power. However, in absence of an expression for the average power, we resort to approximations.

First, we obtain an expression for the average value of the lower bound of the power function of realistic test,
under specific assumptions, 
given in Theorem \ref{Thm:average_power_real_test_h}. 
\begin{theorem}
\label{Thm:average_power_real_test_h}
   In the setting of \textit{Theorem \ref{Thm:lower_bound_power_realistic_h}}, under \textit{Assumptions \ref{Asm:null_uni_dist_h}} and \ref{Asm:alt_uni_dist_h}, when $\epsilon<\min \lbrace a,b-a,1-b \rbrace$,
    \begin{equation}
    \begin{aligned}
\mathbb{E}
\left[
\phi(p')
|p' \in \mathcal{P}_1
\right]
    = 
    \frac{2}{1-(b-a)}
    \left\lbrace
\left(
1-
\frac
{
\epsilon+
\sqrt{\frac{\mathrm{log}r}{r}}
}{b-a} 
\right)^m -
1
    \right\rbrace
    \left(
\epsilon+
\sqrt{\frac{\mathrm{log}r}{r}}
    \right) +
    \left(
1- \frac{2m}{\sqrt{r}}
    \right).
\end{aligned}
\end{equation}
\end{theorem}
Note that the expression involves an unknown $(b-a)$. We thus approximate $(b-a)$ with the difference between the maximum and the minimum values of the estimated Bernoulli parameters, denoted  by $(\hat{b}-\hat{a})$, where $\hat{a}$ and $\hat{b}$ are the outputs of Algorithm \ref{Algo:EstimateRange}. Denoting $H(\epsilon,m,r)
\coloneq
\mathbb{E}[\phi(p')|p' \in \mathcal{P}_1]$,
an approximation for $H(\epsilon,m,r)$ is given by
\begin{equation*}
    \hat{H}(\epsilon,m,r)=
    \frac{2}{1-(\hat{b}-\hat{a})}
    \left\lbrace
\left(
1-
\frac
{
\epsilon+
\sqrt{\frac{\mathrm{log}r}{r}}
}{\hat{b}-\hat{a}} 
\right)^m -
1
    \right\rbrace
    \left(
\epsilon+
\sqrt{\frac{\mathrm{log}r}{r}}
    \right) +
    \left(
1- \frac{2m}{\sqrt{r}}
    \right).
\end{equation*}
for sufficiently large $\tilde{m}$ and $\tilde{r}$. 

We plot the behavior of the lower bound on the average power for the realistic test for $ (a,b) = (0.4, 0.6) $ and $ \alpha = 0.1 $ as a function of $ \epsilon $ for various budgets in the right panel of Figure \ref{fig:analytical-bounds}.
The theoretical results in Section \ref{Sec:Theoretical_results_h} provide justification for Algorithm \ref{Algo:TestOptimal} which maximizes $\hat{H}(\epsilon,m,r)$ with respect to $(\epsilon,m,r)$ satisfying the 
approximate validity constraint
\begin{equation*}
    \left(
1-
\frac{\epsilon-
\sqrt{\frac{\mathrm{log}r}{r}}
}
{\hat{b}-\hat{a}}
    \right)^m + \frac{2m}{\sqrt{r}} \leq \alpha.
\end{equation*}
\begin{corollary}
\label{Cor:asy_valid_consistent_optimal_test}
As budget $\nu \to \infty$ such that $\tilde{m} \to \infty$ and $\tilde{r} \to \infty$, \textit{Algorithm \ref{Algo:TestOptimal}} yields an asymptotically valid sequence of tests. 
\end{corollary}
\section{Experimental Results}
\label{Sec:Numerical_Experiments}
\begin{figure}
    \centering
\includegraphics[width=\linewidth]{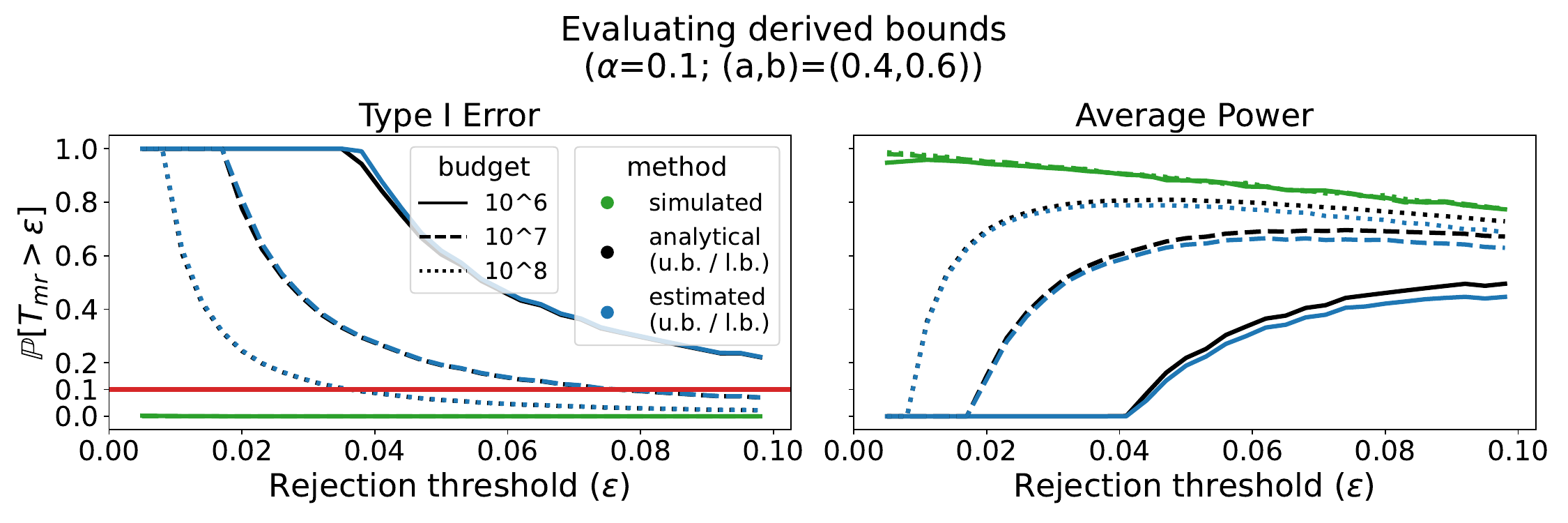}
    \caption{The derived upper bound on Type I Error (left) and the derived lower bound on average power (right) compared to the simulated probability of rejection (green) for different rejection thresholds and various budgets.
    We include both the analytical bound (black) -- where all population parameters are known, and the estimated bound (blue) -- where we use plug-in estimates 
    .
    The derived analytical and estimated bounds properly control both the Type-I Error and Average Power.
    The tightness of both bounds highly depends on the budget.}
    \label{fig:evaluating-bounds}
\end{figure}
We next evaluate the upper bound on the size and the lower bound on the average power of the realistic test and then apply Algorithm \ref{alg:test-optimal} to our motivating example. 
\subsection{Evaluating derived bounds}
\label{Subsec:Simulation}
As in Figure \ref{fig:analytical-bounds}, to evaluate the derived bounds we let $ (a,b) = (0.4, 0.6) $ and $ \alpha = 0.1 $.
We consider $ \epsilon \in \{0.001, 0.002, \hdots, 0.1\} $ and $ \nu \in \{10^{6}, 10^{7}, 10^{8}\} $.
The left panel (respectively right panel) of Figure \ref{fig:evaluating-bounds} includes the derived upper bound on the size (respectively lower bound on the average power) where $ (a,b) $ is known (e.g., the ``analytical" bound) and where $ (a,b) $ must first be estimated (e.g., the ``estimated") bound. 
We also include the simulated probability of rejecting $ H_{0} $ in both panels.
To calculate the simulated probability of rejection we calculate $ m $ according to Eq.~\eqref{eq:m} and set $ r = \nu / m $ for a given $ \epsilon$. 
We sample $ p' \sim \mathrm{Unif}(a,b) $ for size (respectively $p' \sim \mathrm{Unif}((0,a)\cup(b,1)) $ for average power) and then sample $ p_{1}, \hdots, p_{m} $ $i.i.d.$ from $ \mathrm{Unif}(a,b) $.
We estimate each $ p $ with $ r $ samples from Bernoulli($p$), calculate $ T_{m,r} $ per Eq.~\eqref{Eq:realistic_test}, and reject $ H_{0} $ if $ T_{m,r} > \epsilon $.
For a given $ p' $ we repeat the process of sampling $ m $ different $ p_{j} $ from $ \mathrm{Unif}(a,b) $ and estimating each with $ r $ samples from Bernoulli($p_{j}$) 100 times.
The curves labeled ``simulated" in Figure \ref{fig:evaluating-bounds} are the average probabilities of rejecting $ H_{0} $ for $ 1,000 $ different $ p' $.
Finally, the red horizontal line in the left panel corresponds to $ y = \alpha $.

Both panels compel two observations of note: (i) for the budgets under consideration, the estimated bound is close to the analytical bound; that is, using the plug-in estimate of $ (a,b) $ is sufficiently good; and (ii) the derived bounds are relatively loose for all budgets when $ \epsilon $ is small but are tighter for a large budget ($ \nu = 10^7 $ or $10^8$) when $ \epsilon $ is large. 
The closeness of the estimated and analytical bounds is largely due to the fact that we take into account the error in estimating each $ p $ when deriving the bounds and hence rely only on estimation of the interval describing the null region.
The relative looseness of the bounds for small $ \epsilon $ is a general phenomenon-- even for large budgets -- for bounding worst-case scenarios (see, e.g., \cite{alexander-bound}).

\subsection{Revisiting our motivating example}
\begin{figure}
    \centering
\includegraphics[width=\linewidth]{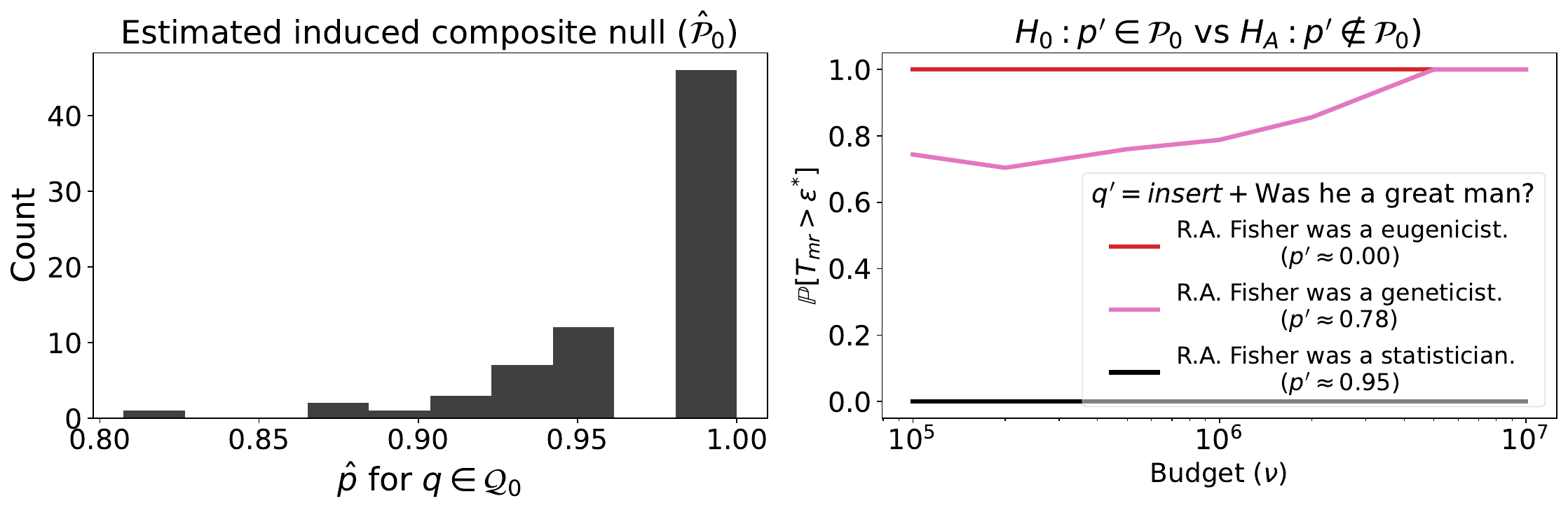}
    \caption{The histogram of estimated Bernoulli parameters of the sampled null queries (left) and the empirical probability of rejecting $ H_{0}: p' \in \mathcal{P}_{0} $ vs. $ H_{A}: p' \notin \mathcal{P}_{0} $ for $ q_{0}~=~\text{``RA Fisher was a statistician. Was he a great man?"}$ using the test described in Alg. \ref{alg:test-optimal} for various $ q' $ (right). 
    The proposed test greatly reduces the undesirable rejections in our motivating example (e.g., when changing ``RA" to ``R.A."), maintains large power for $ p' $ far from $ \mathcal{P}_{0} $ (e.g., when changing ``statistician" to ``eugenicist"), and provides power for $ p' $ close to $ \mathcal{P}_{0} $ when the budget is sufficiently large (e.g., when changing ``statistician" to ``geneticist"). }
    \label{fig:placeholder}
\end{figure}

Consider again the base query $q_{0}=$``RA Fisher was a statistician. Was he a great man?".
As demonstrated in our motivating example, just changing ``RA" to ``R.A." will result in rejecting the equality of the response distributions for large enough $ r $.
To mitigate these types of operationally insignificant rejections via the test described in Eq. \eqref{Eq:realistic_test}, we consider $ \mathcal{Q}_{0} $ to contain the concatenated subelements of the set \{``", ``Prof.", ``Professor"\} $ \times $ \{``RA Fisher", ``R.A. Fisher", ``RA Fisher", `` R.A. Fisher", ``Ronald A Fisher", ``Ronald A. Fisher", ``R A Fisher"\} $ \times $ \{``was a", ``worked as a"\} $ \times $ \{``statistician.", ``biostatistician."\} $ \times $ \{``Was he a great man?"\}; e.g., $ q =$ ``Ronald A. Fisher was a biostatistician. Was he a great man?" $ \in \mathcal{Q}_{0} $.

We let $ f $ be Meta's \texttt{Meta-Llama-3-8B-Instruct} with a temperature $ 1.9 $ and the system prompt ``You are a helpful assistant. You may only respond with `yes' or `no'.".
For each $ q \in \mathcal{Q}_{0} $ we estimate $ \hat{p} $ using $ R = 333,333 $ samples from $ F_{q} $.
The histogram of the estimated elements of $ \mathcal{P}_{0} $ is shown in the left panel of Figure \ref{fig:placeholder}.

We consider three different $ q' $: ``R.A. Fisher was a statistician. Was he a great man?", ``R.A. Fisher was a geneticist. Was he a great man?", and ``R.A. Fisher was a eugenicist. Was he a great man?".
When $ q' =$ ``R.A. Fisher was a statistician. Was he a great man?" we remove it from $ \mathcal{Q}_{0} $.
The other two are different magnitudes of ``farther" from our user-defined notion of of semantic similarity -- changing ``statistician" to ``geneticist" results in a query that is ``closer" to $ \mathcal{Q}_{0}$ than when changing ``statistician" to ``eugenicist".

Following Algorithm \ref{alg:test-optimal},
we sample $ \tilde{m} = 20 $ elements from $ \mathcal{Q}_{0} $ and estimate their corresponding Bernoulli parameters with a random sample of size $ \tilde{r} = 50 $ from the set of $ R $ responses.
We estimate the null region $ (a,b) $ using the unbiased estimates provided in Algorithm \ref{Algo:EstimateRange}.
Then, given a budget $ \nu $ and a level $ \alpha $, we find the optimal triple $ (m^{*}, r^{*}, \epsilon^{*}) $ for $ \epsilon \in \{0.005, 0.01, \hdots, \hat{b} - \hat{a}\} $.
We consider $ \nu \in \{2 \times 10^{4}, 5 \times 10^{4}, 10 \times 10^{4}, 2 \times 10^{5}, \hdots, 10 \times 10^{6}\} $.
For budgets that do not yield a valid test, we use the $ m $ and $ \epsilon $ from the optimal test of the smallest budget that yielded a valid test and reduce $ r $ accordingly.
As an example of an optimal test, for a single instance of the experiment with $ \nu = 5 \times 10^{6} $, $ (\hat{a}, \hat{b}) = (0.898, 1) $ and $ (m^{**}, r^{**}, \epsilon^{**}) = (26, 192307, 0.085) $.

We report the average probability of rejecting $ H_{0}: p' \in \mathcal{P}_{0} $ in the right panel of Figure \ref{fig:placeholder} as a function of budget.
The average is over $ 250 $ different instances of the experiment -- e.g., estimating $ (\hat{a}, \hat{b}) $, finding the optimal tests, sampling $ r^{**} $ from the set of  $ 333,333 $ responses for query sampled query  -- except for when $ \nu = 10 \times 10^{6} $ because $ r^{**} \approx 333,333 $.
For $ \nu = 10 \times 10^{6} $, the reported average probability of rejection is the average over tests resulting from different $ \tilde{m} $ and $ \tilde{r} $ corresponding optimal tests.

Our proposed test properly controls the operationally insignificant rejections described in the introduction and maintains non-trivial power both when changing ``statistician" to ``geneticist" and ``statistician" to ``eugenicist" even when $ \nu $ is small. 
Notably, the proposed test has more power for the query that is farther from $ \mathcal{Q}_{0} $ than the query closer to $ \mathcal{Q}_{0} $.
We also note that the smallest budget where the lower bound on the size of the test is less than $ \alpha = 0.1 $ is $ \nu = 2 \times 10^{6} $ -- which likely causes the power to be outsized when the budget is particularly small.
The increase in power likely comes with an increase in size when considering a larger $ \mathcal{Q}_{0} $, though we do not observe it here.

\subsection{Cost}
The above experiment was conducted by prompting Meta's \texttt{Meta-Llama-3-8B-Instruct} approximately $ 72 \cdot 333,333 \approx 2 \times 10^{7} $ times.
Generating all the responses took approximately $ 100 $ hours on a single Nvidia H100 and cost approximately $ \$200 $.
We discuss extensions of our proposed test -- such as more intelligent sampling of $ q \in \mathcal{Q}_{0} $ or considering different $ r $ for different $ q $ -- in the discussion below.
\section{Discussion}
\label{Sec:Discussion}

In this paper, we introduce a statistical framework for testing the difference of response distributions in the context of semantically irrelevant perturbations of a base query. 
We restrict ourselves to the regime of
responses. Motivated by \cite{bodmer2021outstanding}'s discussion about the famous scientist Ronald A. Fisher's contribution to statistics and genetics alongside his controversial views on eugenics,
we use our methodology to test the differences in distributions of LLM responses to queries pertaining to Ronald A. Fisher's
identity as statistician and eugenicist. 
Our investigation of statistical hypothesis tests deployed to detect significant changes in the response distribution due to query perturbation is in line with the principle of stability of statistical results to reasonable perturbations in the data, discussed in \cite{agarwal2025pcs} and \cite{yu2024veridical}.

Recall that we reject $H_0:p' \in \mathcal{P}_0$ when the test statistic $T_{m,r}$ is larger than a threshold $\epsilon$, and the quantities $\epsilon$, $m$ and $r$ are chosen such that they satisfy a (approximate) validity constraint. The expression for the validity constraint is deduced based on the assumption that the Bernoulli parameters of the sampled null queries are distributed uniformly.  
In reality, the distribution of the Bernoulli parameters of the sampled queries is unknown, because the map from the query set $\mathcal{Q}$ to the set of corresponding Bernoulli parameters, $\mathcal{P}$, is unknown.   
The histogram of the estimated Bernoulli parameters $\hat{p}_j$
in Figure \ref{fig:placeholder} indicate non-uniformity in the distribution of the Bernoulli parameters,
providing motivation for future work generalizing our setting to, for example, the case
where the Bernoulli parameters follow a mixture of Beta distributions.

We deal herein with a black-box setting because in reality, information about the internal structure of the model is often not available. However, our ideal test can be used to deal with a white-box setting where the user has access to the internal structure of the LLM.

Since $\mathcal{P}_0$ is unknown, we repeatedly draw samples from it to ensure sufficient coverage of $\mathcal{P}_0$. 
We assumed that the null queries are being sampled independently.
Developing methods for intelligent sampling to ensure adequate coverage at a lower cost may be a promising direction for future work.

Our proposed test is based on the idea that the realistic test statistic $T_{m,r}$ approximates the ideal test statistic $\Tilde{T}_m$ for sufficiently large $r$, and that in order to have high power for the ideal test, $\epsilon$ must be small, which warrants large $m$. Thus, for a given level of significance, in order to have asymptotic validity and consistency, it is necessary to have $r=\omega(m^2)$ while $m \to \infty$, which is established in Corollary \ref{Cor:asy_valid_const_true_par_h}. 

It is perhaps reasonable to assume some Lipschitz-like continuity property
for the LLM map from query to response distribution.
Letting $c_0$ be the local Lipschitz constant associated with the base query $q_0$
and $d_\mathcal{Q}$ be some distance on query strings (e.g., Levenshtein distance,
or a bespoke distance capturing user-defined semantically irrelevant query perturbations),
this suggests
$$|p-p_0| \leq c_0 d_\mathcal{Q}(q,q_0).$$
Generalizing from Bernoullis to arbitrary response distributions $F$ equipped with some appropriate distance $d$ (such as total variation),
this becomes $d(F,F_0) \leq c_0 d_\mathcal{Q}(q,q_0)$.
In the Bernoulli case addressed in this paper,
if $$\mathcal{Q}_0 \stackrel{\text{def}}{=} \{q: d_\mathcal{Q}(q,q_0) \leq \epsilon\}$$
then $$\mathcal{P}_0 \subset \{p: |p-p_0| \leq c_0 \epsilon\}$$
and thus we could consider 
$H_0^{Lip}: p \in p_0 \pm c_0 \epsilon$ --
a valid test for $H_0^{Lip}$ is valid for our original $H_0: p \in \mathcal{P}_0$.
With $c_0$ known, a straightforward variation of classical two-sample Neyman-Pearson testing for equality of two Bernoulli parameters applies for $H_0^{Lip}$,
providing a simple and compelling illustration the utility of modeling LLM maps as Lipschitz.
However, with $c_0$ unknown this formulation must contend with precisely the same complication that we have addressed in this paper --
an {\em unknown} range for the null probabilities induced by $\mathcal{Q}_0$.

\begin{figure}[h!]
    \centering
    \includegraphics[scale=0.30]{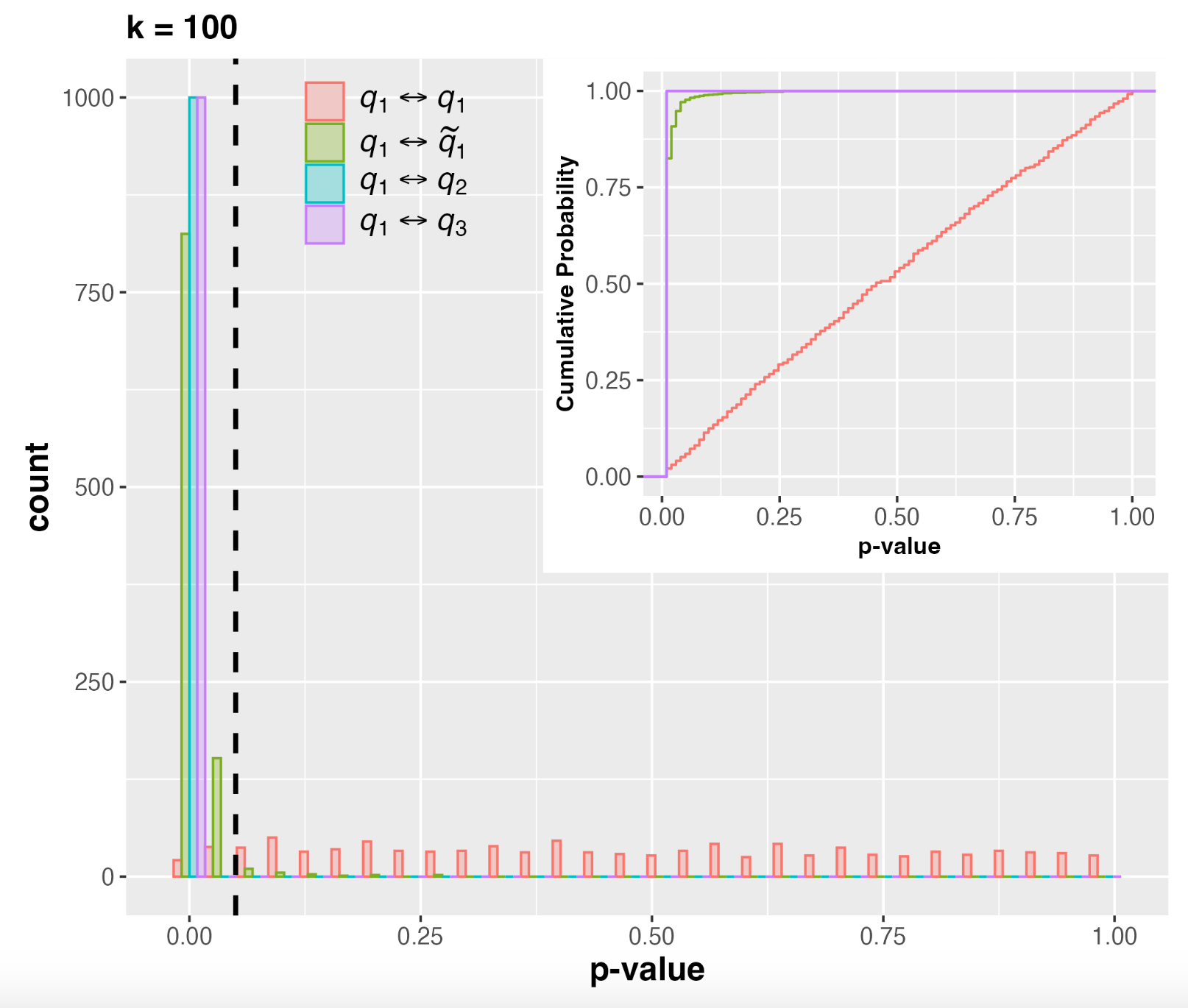}
    \caption{Distribution of p-values for tests for semantic equivalence of queries, in the setting of general (non-binary) responses. The large language model used is \texttt{google/gemma/2-2b-it} and the embedding function $g$ is \texttt{nomic-ai/nomic-embed-text-v2-moe}.
    For every query,  we bootstrap $k=100$ responses from a pool of $1000$ randomly generated responses, and implement Szekely's Energy Test on any pair of queries, obtaining a p-value. We repeat this procedure on $m=1000$ Monte Carlo samples to obtain an empirical distribution of p-values, which is shown in the figure.}
    \label{fig:Energy_general}
\end{figure}

Indeed, a scope for future extension of our work involves the investigation of a regime of generalized responses, instead of a regime of binary responses. In such a generalized regime, if the vectorized versions of the responses can be modeled with parametric distributions, then a path similar to ours can be followed. However, in a distribution-free setting, a radically different approach will be needed. We demonstrate an approach for testing exact semantic equivalence between two queries with a simple null hypothesis, and hope for future extension along this line for testing semantic similarity involving a composite null hypothesis. We use Szekely's energy test \citep{Szkely2004TESTINGFE} to test $H_0:F_{q_1}=F_{q_2}$ where $F_{q_1}$ and $F_{q_2}$ are respectively the response distributions induced by the queries $q_1$ and $q_2$. In our example, we use the following queries: 
\newline
\hspace*{0.25in}$q_1=\text{``Describe why you think RA Fisher was a great statistician''}$;
\newline
\hspace*{0.25in}$\Tilde{q}_1=\text{``Describe why you think R.A. Fisher was a great statistician''}$;
\newline
\hspace*{0.25in}$q_2=\text{``Describe why you think RA Fisher was a great geneticist''}$;
\newline
\hspace*{0.25in}$q_3=\text{``Describe why you think RA Fisher was a great eugenicist''}$. 
\newline
For every pair of queries, we plot the empirical distribution of p-values in Figure \ref{fig:Energy_general}.
We find that for testing $H_0:F_{q_1}=F_{\Tilde{q}_1}$ where the perturbation is to be considered semantically irrelevant, the proportion of rejection of $H_0$ is undesirably high. This shows the need for further investigation into the regime of generalized responses
for the case of a composite null consisting of
semantically irrelevant query perturbations.
\newline

Finally, any analysis in this paper can be applicable to generative models in general, including but not restricted to large language models.

\newpage
\acks{Support for this effort provided by Defense Advanced Research Projects Agency (DARPA) Artificial Intelligence Quantified (AIQ) award number HR00112520026. We would also like to thank Robert O. Ness and Youngser Park for their help on the numerical experiments.
}


\clearpage

\bibliography{reff}

\section{Appendix: Proofs of results}
\label{sec:Proofs}
\textit{
\textbf{Lemma 1.} 
Suppose that for every query $ q \in \{q_{1}, \hdots, q_{m}, q'\} $ we observe iid replicates of responses denoted by
$f(q)_1,\dots,
f(q)_r \sim^{iid} \mathrm{Bernoulli}(p)$  where $p$ is the Bernoulli parameter of the query $q$. 
Define 
\begin{equation*}
\begin{aligned}
    &\hat{p}_j=\frac{1}{r} \sum_{k=1}^r f(q_j)_k \text{ for all $j \in [m]$}, \\
    &\hat{p}'=\frac{1}{r} \sum_{k=1}^r f(q')_k, 
\end{aligned}
\end{equation*}
and set $\tilde{T}_m=\min_{j \in [m]}|p_j-p'|$, $T_{m,r}=\min_{j \in [m]} |\hat{p}_j-\hat{p}'|$. 
Then, conditioning on $ (p_{1},\dots p_{m}, p')$,
\begin{equation*}
\mathbb{P}
\left[
   \big|
T_{m,r}-\tilde{T}_m
   \big|
   <
   \sqrt{\frac{\mathrm{log}r}{r}}
   \bigg|
   (p_1,\dots p_{m},p')
\right] \geq 1 -
\frac{2m}{\sqrt{r}}
.
\end{equation*}
}
\newline
\textbf{Proof.} 
Note that,
\begin{equation*}
\begin{aligned}
\bigg|
    \big|
\frac{1}{r} \sum_{k=1}^r X^{(j)}_k - 
\frac{1}{r} \sum_{k=1}^r X'_k
    \big| -
    \big|
p_j-p'
    \big|
\bigg| 
&\leq  
\bigg|
\left(
\frac{1}{r} \sum_{k=1}^r X^{(j)}_k-
\frac{1}{r} \sum_{k=1}^r X'_k
\right)-
(p_j-p') 
\bigg| \\
&=
\bigg|
\sum_{k=1}^r 
\frac{
X^{(j)}_k-X^{'}_k
}{r}-
\mathbb{E}
\left(
\sum_{k=1}^r 
\frac{
X^{(j)}_k-X^{'}_k
}{r}
\right)
\bigg|.
\end{aligned}
\end{equation*}
Now, note that for every $j \in [m]$, for all $k \in [r]$,
$
\mathbb{P}
\left[
\frac{
X^{(j)}_k-X^{'}_k
}{r} \in [-\frac{1}{r},\frac{1}{r}]
\right]=1
$. Thus, using \textit{Lemma 1}, 
\begin{equation*}
\begin{aligned}
&\mathbb{P}
\left[
\bigg|
\sum_{k=1}^r 
\frac{
X^{(j)}_k-X^{'}_k
}{r}-
\mathbb{E}
\left(
\sum_{k=1}^r 
\frac{
X^{(j)}_k-X^{'}_k
}{r}
\right)
\bigg|
< t
\right]
\geq 1 - 2e^{
-\frac{r t^2}{2}
} \\
&\implies 
\mathbb{P}
\left[
\bigg|
    \big|
\frac{1}{r} \sum_{k=1}^r X^{(j)}_k - 
\frac{1}{r} \sum_{k=1}^r X'_k
    \big| -
    \big|
p_j-p'
    \big|
\bigg| 
< t
\right] \geq 
1 - 2e^{
-\frac{r t^2}{2}
}
\end{aligned}
\end{equation*}
\newline
Choosing $t=\sqrt{\frac{\mathrm{log}r}{r}}$, we get for all $j \in [m]$,
\begin{equation*}
    \mathbb{P}
    \left[
\bigg|
    \big|
\frac{1}{r} \sum_{k=1}^r X^{(j)}_k - 
\frac{1}{r} \sum_{k=1}^r X'_k
    \big| -
    \big|
p_j-p'
    \big|
\bigg| 
< \sqrt{\frac{\mathrm{log}r}{r}}
    \right]
\geq 1-2r^{-\frac{1}{2}}.
\end{equation*}
We know that 
$
|\min_x f(x)-\min_x g(x)|
\leq \max_x |f(x)-g(x)|
$. Thus,
\begin{equation*}
\begin{aligned}
    |T_{m,r}-\tilde{T}_m| 
    &=
    \bigg|
    \min_{j \in [m]}
    \big|
\frac{1}{r} \sum_{k=1}^r X^{(j)}_k - 
\frac{1}{r} \sum_{k=1}^r X'_k
    \big| -
    \min_{j \in [m]}
    \big|
p_j-p'
    \big|
\bigg|  \\
&\leq 
\max_{j \in [m]}
\bigg|
    \big|
\frac{1}{r} \sum_{k=1}^r X^{(j)}_k - 
\frac{1}{r} \sum_{k=1}^r X'_k
    \big| -
    \big|
p_j-p'
    \big|
\bigg|. 
\end{aligned}
\end{equation*}
Now, 
\begin{equation*}
\begin{aligned}
    &\mathbb{P}
    \left[
\left|    T_{m,r}-\tilde{T}_m
\right|< \sqrt{\frac{\mathrm{log}r}{r}}
\bigg| (p_1,\dots,p_m,p')
    \right] \\ &\geq 
    \mathbb{P}
    \left[
    \max_{j \in [m]}
\bigg|
    \big|
\frac{1}{r} \sum_{k=1}^r X^{(j)}_k - 
\frac{1}{r} \sum_{k=1}^r X'_k
    \big| -
    \big|
p_j-p'
    \big|
\bigg| <
\sqrt{\frac{\mathrm{log}r}{r}}
\bigg|(p_1,\dots,p_m,p')
    \right] \\
    &=
    1-
    \mathbb{P}
    \left[
    \max_{j \in [m]}
\bigg|
    \big|
\frac{1}{r} \sum_{k=1}^r X^{(j)}_k - 
\frac{1}{r} \sum_{k=1}^r X'_k
    \big| -
    \big|
p_j-p'
    \big|
\bigg| \geq 
\sqrt{\frac{\mathrm{log}r}{r}}
\bigg| (p_1,\dots,p_m,p')
    \right] \\
    &\geq 1- 
    \sum_{j=1}^m 
    \mathbb{P}
    \left[
\bigg|
    \big|
\frac{1}{r} \sum_{k=1}^r X^{(j)}_k - 
\frac{1}{r} \sum_{k=1}^r X'_k
    \big| -
    \big|
p_j-p'
    \big|
\bigg| \geq \sqrt{\frac{\mathrm{log}r}{r}}
\bigg| (p_1,\dots,p_m,p')
    \right] \\
    &\geq 1-\sum_{j=1}^m 2 r^{-\frac{1}{2}} \\
    &\geq 1 -\frac{2m}{\sqrt{r}}.
\end{aligned}
\end{equation*}
\newline
\newline
\textit{
\textbf{Lemma 2.}
Under \textit{Assumption \ref{Asm:null_alt_independence_h}} and \textit{Assumption \ref{Asm:null_uni_dist_h}}  
on $ q_{1}, \hdots, q_{m} $, suppose we observe the true corresponding Bernoulli parameters $ p_{1}, \hdots, p_{m}, p' $ and we reject $H_0$ if $\tilde{T}_m>\epsilon$.
Then, if $\epsilon<\min \lbrace a,b-a,1-b \rbrace$, a sufficient condition to ensure that the test is valid at level of significance $\alpha$ is given by
\begin{equation}
\label{eq:m}
    m \geq 
    \Bigl\lceil
    \frac
    {
|\mathrm{log}(\alpha)|
    }
    {
    |\mathrm{log}(1-\frac{\epsilon}{b-a})|
    }
    \Bigl\rceil
    .
\end{equation}
}
\newline
\textbf{Proof.} If the true $p_j$ and $p'$ were available, we should  choose $m$ such that 
\newline
$\sup_{p' \in \mathcal{P}_0} \mathbb{P}[\tilde{T}_m>\epsilon|p'] \leq \alpha$.
\newline
We divide all possibilities into three cases viz. $(b-a)<\min \lbrace a,1-b \rbrace$, $a<(b-a)<(1-b)$ and $(1-b)<(b-a)<a$. 
\newline
\newline
\underline{Case I:}($b-a<\min \lbrace a,1-b \rbrace$)
\newline
\newline
Here, $b-a=\min \lbrace a,b-a,1-b \rbrace$. 
\newline
For $p' \in \mathcal{P}_0=[a,b]$,
when $\epsilon<\frac{b-a}{2}$,
\begin{equation*}
\begin{aligned}
    \mathbb{P}
    \left[
    \tilde{T}_m>\epsilon|p'
    \right] &=
    \mathbb{P}
    \left[
    \min_{j \in [m]}
    |p_j-p'|>\epsilon \bigg|
    p'
    \right] \\
    &= \prod_{j=1}^m \mathbb{P}
    \left[ 
    |p_j-p'|>\epsilon
    \bigg| p'
    \right] \\
    &= 
    \prod_{j=1}^m
    \left(
    1-\mathbb{P}
    \left[
    |p_j-p'|
    \leq \epsilon
    \bigg| p'
    \right]
    \right) \\
    &=
       \left\{
\begin{array}{ll}
      \left(
      1-\frac{2 \epsilon}{b-a}
      \right)^m,
       & p' \in (a+\epsilon,b-\epsilon) \\
       \left(
      \frac{b-p'-\epsilon}{b-a}
      \right)^m,
       & p' \in (a,a+\epsilon) \\
      \left(\frac{p'-\epsilon-a}{b-a}
      \right)^m, & p' \in (b-\epsilon,b) \\
\end{array} 
\right. \\
&\leq 
 \left\{
\begin{array}{ll}
      \left(
      1-\frac{2 \epsilon}{b-a}
      \right)^m,
       & p' \in (a+\epsilon,b-\epsilon) \\
       \left(1-
\frac{\epsilon}{b-a}
      \right)^m,
       & p' \in (a,a+\epsilon) \\
      \left(1-\frac{\epsilon}{b-a}
      \right)^m, & p' \in (b-\epsilon,b) \\
\end{array} 
\right. 
\end{aligned}
\end{equation*}
For $p' \in \mathcal{P}_0$, when 
$\frac{b-a}{2}\leq \epsilon<(b-a)$, we have,
\begin{equation*}
\begin{aligned}
\vspace{-3cm} 
    \mathbb{P}
    \left[
\tilde{T}_m>\epsilon|p'
    \right] &=
    \mathbb{P}
    \left[
    \min_{j \in [m]}
    |p_j-p'|>\epsilon |p'
    \right] \\
    &= \prod_{j=1}^m 
    \mathbb{P}
    \left[
|p_j-p'|>\epsilon \bigg|p'
    \right] \\
    &= \prod_{j=1}^m 
    \left(
    1-
    \mathbb{P}
    \left[
    |p_j-p'| \leq \epsilon
    \bigg| p'
    \right]
    \right) \\
    &= \prod_{j=1}^m 
    \left(
    1-
\int_{p'-\epsilon}^{p'+\epsilon}
f_j(x)dx
    \right) \\
    &=
    \left\{
\begin{array}{ll}
   \left(
      \frac{b-p'-\epsilon}{b-a}
      \right)^m,
       & p' \in (a,b-\epsilon] \\
       0,
       & p' \in (b-\epsilon,a+\epsilon) \\
      \left(\frac{p'-\epsilon-a}{b-a}
      \right)^m, & p' \in [a+\epsilon,b)
\end{array}
    \right. \\
&\leq 
\left\{
\begin{array}{ll}
 \left(1-
      \frac{\epsilon}{b-a}
      \right)^m,
       & p' \in (a,b-\epsilon] \\
       0,
       & p' \in (b-\epsilon,a+\epsilon) \\
      \left(1-\frac{\epsilon}{b-a}
      \right)^m & p' \in [a+\epsilon,b).
\end{array}
\right.
\end{aligned}
\end{equation*}
Thus, when $\epsilon<(b-a)=\min \lbrace a,b-a,1-b \rbrace$, for all $p' \in \mathcal{P}_0$,
$
\mathbb{P}
    \left[
  \tilde{T}_m>\epsilon|p'  
    \right] \leq 
    \left(
1-\frac{\epsilon}{b-a}
    \right)^m
$.
\newline
\newline
\underline{Case II:}($a<(b-a)<(1-b)$)
\newline
\newline
Here, $a=\min \lbrace a,b-a,1-b \rbrace$. 
\newline
When $\epsilon<a$, for $p' \in \mathcal{P}_0$, 
\begin{equation*}
\begin{aligned}
    \mathbb{P}
    \left[
    \tilde{T}_m>\epsilon|p'
    \right] &= 
    \mathbb{P}
    \left[
    \min_{j \in [m]}
    |p_j-p'|>\epsilon
    \bigg|p'
    \right] \\
    &= \prod_{j=1}^m 
    \mathbb{P}
    \left[
    |p_j-p'|>\epsilon
    \bigg|p'
    \right] \\
    &= 
    \prod_{j=1}^m 
    \left(
    1-\mathbb{P}
    \left[
    |p_j-p'| \leq \epsilon \bigg| p'
    \right]
    \right) \\
    &= \prod_{j=1}^m 
    \left(1- 
\int_{p'-\epsilon}^{p'+\epsilon}
f_j(x)dx\right) \\
&= 
\left\{
\begin{array}{ll}
\left(1-\int_{a}^{p'+\epsilon}
\frac{1}{b-a}dx \right)^m,
& p' \in (a,a+\epsilon] \\
\left(1- \int_{p'-\epsilon}^{p'+\epsilon}
\frac{1}{b-a}dx \right)^m,
& p' \in (a+\epsilon,b-\epsilon) \\
\left(1-\int_{p'-\epsilon}^{b} \frac{1}{b-a} \right)^m, 
& p' \in [b-\epsilon,b) \\
\end{array}
\right. \\
&= \left\{
\begin{array}{ll}
\left(
\frac{b-p'-\epsilon}{b-a} \right)^m,
& p' \in (a,a+\epsilon] \\
\left(1- \frac{2 \epsilon}{b-a}
\right)^m, & p' \in (b-\epsilon,a+\epsilon) \\
\left(\frac{p'-\epsilon-a}{b-a}
\right)^m, & p' \in [a+\epsilon,b) \\
\end{array}
\right. \\
&\leq 
\left\{
\begin{array}{ll}
\left(1-\frac{\epsilon}{b-a}\right)^m,
& p' \in (a,a+\epsilon] \\
\left(1-\frac{2 \epsilon}{b-a}
\right)^m,
& p' \in (b-\epsilon,a+\epsilon) \\
\left(1-\frac{\epsilon}{b-a}
\right)^m, & p' \in [a+\epsilon,b) \\
\end{array}
\right.
\end{aligned}
\end{equation*}
Thus, when $\epsilon<a=\min \lbrace a,b-a,1-b \rbrace$, for all $p' \in \mathcal{P}_0$, we have 
$
\mathbb{P}
\left[
\tilde{T}_m>\epsilon|p'
\right] \leq 
\left(
1- \frac{\epsilon}{b-a}
\right)^m
$.
\newline
\underline{Case III:}($1-b<b-a<a$)
\newline
\newline
Here, $1-b= \min \lbrace a,b-a,1-b \rbrace$. 
\newline
For $p' \in \mathcal{P}_0=[a,b]$,
when  $\epsilon<(1-b)$,
\begin{equation*}
\begin{aligned}
    \mathbb{P}
    \left[
\tilde{T}_m>\epsilon|p'
    \right] &=
    \prod_{j=1}^m 
    \left(
1- \int_{p'-\epsilon}^{p'+\epsilon
} f_j(x)dx
    \right) \\
    &=
    \left\{
\begin{array}{ll}
      \left(1-
      \int_{a}^{p'+\epsilon}
      \frac{1}{b-a}dx
      \right)^m,
       & p' \in (a,a+\epsilon] \\
       \left(1-
       \int_{p'-\epsilon}^{p'+\epsilon}
       \frac{1}{b-a}dx
       \right)^m,
       & p' \in (a+\epsilon,b-\epsilon) \\
      \left(1-
      \int_{p'-\epsilon}^b
      \frac{1}{b-a}dx
      \right)^m, & p' \in [b-\epsilon,b) \\
\end{array}
\right. \\
&=
\left\{
\begin{array}{ll}
      \left(
      \frac{b-p'-\epsilon}{b-a}
      \right)^m,
       & p' \in (a,a+\epsilon] \\
       \left(1-
       \frac{2 \epsilon}{b-a}
       \right)^m,
       & p' \in (a+\epsilon,b-\epsilon) \\
      \left(
      \frac{p'-\epsilon-a}{b-a}
      \right)^m, & p' \in [b-\epsilon,b) \\
\end{array}
\right. \\
&\leq 
\left\{
\begin{array}{ll}
      \left(1-
      \frac{\epsilon}{b-a}
      \right)^m,
       & p' \in (a,a+\epsilon] \\
       \left(1-
       \frac{2 \epsilon}{b-a}
       \right)^m,
       & p' \in (a+\epsilon,b-\epsilon) \\
      \left(1-
      \frac{\epsilon}{b-a}
      \right)^m, & p' \in [b-\epsilon,b) \\
\end{array}
\right. 
\end{aligned}
\end{equation*}
Thus, when $\epsilon<b=\min \lbrace a,b-a,1-b \rbrace$, for all $p' \in \mathcal{P}_0$,
\begin{equation*}
    \mathbb{P}
    \left[
\tilde{T}_m>\epsilon|p'
    \right]
    \leq 
    \left(
1-\frac{\epsilon}{b-a}
    \right)^m.
\end{equation*}
\newline
Combining all three cases, when $\epsilon< \min \lbrace a,b-a,1-b \rbrace$, for all $p' \in \mathcal{P}_0$,
\begin{equation*}
    \mathbb{P}
    \left[
    \tilde{T}_m>\epsilon|p'
    \right] \leq 
    \left(
1- \frac{\epsilon}{b-a}
    \right)^m. 
\end{equation*}
\newline
Hence, to have $\sup_{p' \in \mathcal{P}_0} \mathbb{P}[\tilde{T}_m>\epsilon|p'] \leq \alpha$, it suffices to ensure that
\begin{equation*}
    m \geq 
    \frac
    {
|\mathrm{log}(\alpha)|
    }
    {
    |\mathrm{log}(1-\frac{\epsilon}{b-a})|
}.
\end{equation*}
\textit{
\textbf{Lemma 3.} 
In our setting, under \textit{Assumptions} 
\ref{Asm:null_alt_independence_h} and
\ref{Asm:null_uni_dist_h}, suppose we observe the true Bernoulli parameters, and we want to test $H_0:p' \in \mathcal{P}_0$ versus $p' \notin \mathcal{P}_0$ at level of significance $\alpha$. Our ideal decision rule rejects $H_0$ if $\tilde{T}_m>\epsilon$ for some chosen threshold $\epsilon$, where $\tilde{T}_m=\min_{j \in [m]} |p_j-p'|$. 
For $\epsilon<\min \lbrace a,b-a,1-b \rbrace$, the power function of the ideal test is given by
\begin{equation*}
\begin{aligned}
\tilde{\beta}(p')
\equiv \tilde{\beta}_{m,\epsilon}(p')
\coloneq
    \mathbb{P}
    \left[
    \tilde{T}_m>\epsilon
    \bigg| p'
    \right]=
    \left\{
\begin{array}{ll}
      1,
       & p'\in (0,a-\epsilon] \cup [b+\epsilon,1) \\
       \left(
      \frac{b-p'-\epsilon}{b-a}
      \right)^m,
       & p' \in (a-\epsilon,a) \\
       \left(
      \frac{p'-\epsilon-a}{b-a}
      \right)^m, & p' \in (b,b+\epsilon)  \\
\end{array} 
\right.
\end{aligned}
\end{equation*}
}
\newline
\textbf{Proof.} 
For $\epsilon<\min \lbrace a,b-a,1-b \rbrace$,
we can see,
\begin{equation*}
\begin{aligned}
    \mathbb{P}
    \left[
    \tilde{T}_m>\epsilon
    \bigg| p'
    \right] &=
    \mathbb{P}
    \left[
    \min_{j \in [m]}
    |p_j-p'|>\epsilon
    \bigg|p'
    \right]  \\
    &= 
    \prod_{j=1}^m
    \mathbb{P}
    \left[
    |p_j-p'|>\epsilon
    \bigg|p' 
    \right] \\
    &= \prod_{j=1}^m
    \bigg(
    1-
    \mathbb{P}
    \left[
    |p_j-p'| \leq \epsilon
    \bigg| p'
    \right]
    \bigg) \\
    &= \prod_{j=1}^m 
    \bigg(
1-\int_{p'-\epsilon}^{p'+\epsilon}
f_j(x)dx
    \bigg) \\
    &=
    \left\{
\begin{array}{ll}
      1,
       & p'\in (0,a-\epsilon] \cup [b+\epsilon,1) \\
       \left(
      \frac{b-p'-\epsilon}{b-a}
      \right)^m,
       & p' \in (a-\epsilon,a) \\
       \left(
      \frac{p'-\epsilon-a}{b-a}
      \right)^m, & p' \in (b,b+\epsilon)  \\
\end{array} 
\right.
    \end{aligned}
\end{equation*}
\newline
\newline
\textit{
\textbf{Theorem 1.} Suppose \textit{Assumption \ref{Asm:null_alt_independence_h}} holds, and consider 
the setting of \textit{Lemma \ref{lemma:concentration_of_test_statistic}}. For any $p' \in (0,1)$,
\begin{equation*}
    \mathbb{P}
    \left[
    \left|
T_{m,r}-\tilde{T}_m
    \right|
    < 
    \sqrt{\frac{\mathrm{log}r}{r}}
    \bigg| p'
    \right]
    \geq 
1- \frac{2m}{\sqrt{r}}.
\end{equation*}
}
\textbf{Proof.} 
It is easy to see,
\begin{equation*}
\begin{aligned}
    &\mathbb{P}
\left[\left|
T_{m,r}-\tilde{T}_m
\right| < 
\sqrt{\frac{\mathrm{log}r}{r}}
\bigg| p' \right] \\
&= \int_{p_1=a}^b \int_{p_2=a}^b \dots \int_{p_m=a}^b
\mathbb{P}
\left[ \left|
T_{m,r}-\tilde{T}_m \right| < 
\sqrt{\frac{\mathrm{log}r}{r}}
\bigg|(p_1,\dots p_m,p')
\right]
f_1(p_1)f_2(p_2)\dots f_m(p_m) dp_1 \dots dp_m \\
&=
1- \frac{2m}{\sqrt{r}}.
\end{aligned}
\end{equation*}
\textit{
\textbf{Theorem 2.} Suppose 
\textit{Assumption \ref{Asm:null_alt_independence_h}} and
\textit{Assumption \ref{Asm:null_uni_dist_h}} hold, and the Bernoulli parameters $p_1,\dots,p_m$ and $p'$ are not observed. For every query $q \in \lbrace
q_1,\dots,q_m,q'
\rbrace$, iid replicates of responses, denoted by $f(q)_1,\dots,f(q)_r$ are obtained. Define $\hat{p}_j=\frac{1}{r} \sum_{k=1}^r f(q_j)_k$ and $\hat{p}'=\frac{1}{r} \sum_{k=1}^r f(q')_k$. 
Recall that
for testing $H_0:p' \in \mathcal{P}_0$ versus $H_A:p' \notin \mathcal{P}_0$,
our decision rule rejects $H_0$ if $T_{m,r}>\epsilon$, for a chosen threshold $\epsilon$, where 
$T_{m,r}=\min_{j \in [m]}|\hat{p}_j-\hat{p}'|$. 
 When $\epsilon<\min \lbrace a,b-a,1-b \rbrace$ and
$r$ is sufficiently large, for all $p' \in \mathcal{P}_0$, 
\begin{equation*}
\begin{aligned}
    \mathbb{P}
    \left[
    T_{m,r}>\epsilon
    \bigg| p'
    \right] 
    \leq 
    \left(
1- 
\frac
{
\epsilon -
\sqrt{\frac{\mathrm{log}r}{r}}
}{b-a}
    \right)^m +
    \frac{2m}{\sqrt{r}}.
\end{aligned}
\end{equation*}
}
\newline
\textbf{Proof.}
Observe that, when $\epsilon< \min \lbrace a,b-a,1-b \rbrace$ and $r$ is sufficiently large, for all $p' \in \mathcal{P}_0$,
we have
\vspace{-3cm}
\begin{equation*}
\begin{aligned}
    \mathbb{P}
    \left[
    T_{m,r}>\epsilon
    \bigg| p'
    \right]
    &=
    \mathbb{P}
    \left[
    T_{m,r}>\epsilon,
    \big|
T_{m,r}-\tilde{T}_m
    \big|
    <\sqrt{\frac{\mathrm{log}r}{r}}
    \bigg| p'
    \right] +
    \mathbb{P}
    \left[
    T_{m,r}>\epsilon,
    \big|
T_{m,r}-\tilde{T}_m
    \big|
    \geq
\sqrt{\frac{\mathrm{log}r}{r}}
    \bigg| p'
    \right] \\
    &\leq 
    \mathbb{P}
    \left[
    \tilde{T}_m>
    \epsilon-
    \sqrt{\frac{\mathrm{log}r}{r}} \bigg|p'
    \right] +
    \mathbb{P}
    \left[
    \left|
T_{m,r}-\tilde{T}_m
\right|
\geq 
\sqrt{\frac{\mathrm{log}r}{r}} 
    \bigg| p'
    \right] \\
    &\leq 
    \left(
1- 
\frac
{
\epsilon- \sqrt{\frac{\mathrm{log}r}{r}}
}{b-a}
    \right)^m +
    \frac{2m}{\sqrt{r}}
    \hspace{01cm}
    \text{[using Theorem \ref{Thm:lower_bound_on_real_size_h} and Theorem \ref{thm:concentration_of_test_statistic_only_ew_conditioning}]}.
\end{aligned}
\end{equation*}
\vspace{-3cm}
\newline
\textit{
\textbf{Theorem 3.} Consider the setting of \textit{Theorem \ref{Thm:lower_bound_on_real_size_h}}. 
Recall that for testing $H_0:p' \in \mathcal{P}_0$ versus $H_1:p' \notin \mathcal{P}_0$, our realistic decision rule rejects $H_0$ if $T_{m,r}>\epsilon$, for a chosen threshold $\epsilon$, where 
$T_{m,r}=\min_{j \in [m]}|\hat{p}_j-\hat{p}'|$.  When $\epsilon<\min \lbrace a,b-a,1-b \rbrace$ and $r$ is sufficiently large,
$\mathbb{P}[T_{m,r}>\epsilon|p'] \geq \phi(p')$ where the lower bound $\phi$ is given by
\begin{equation*}
\begin{aligned}
\phi(p') \coloneq
\left\{
\begin{array}{ll}
\left( 1-\frac{2m}{\sqrt{r}}
\right),
& p'\in \left(0,a-\epsilon-
\sqrt{\frac{\mathrm{log}r}{r}}
\right]  \cup 
\left[b+\epsilon+\sqrt{\frac{\mathrm{log}r}{r}},1 \right) \\
\left(1-\frac{\epsilon+
\sqrt{\frac{\mathrm{log}r}{r}}
}{b-a}\right)^m
- \frac{2m}{\sqrt{r}}, & p' \in (a-\epsilon-
\sqrt{\frac{\mathrm{log}r}{r}},a) \\
\left(1-\frac{\epsilon+
\sqrt{\frac{\mathrm{log}r}{r}}
}{b-a}\right)^m
- \frac{2m}{\sqrt{r}}, 
& p' \in (b,b+\epsilon+
\sqrt{\frac{\mathrm{log}r}{r}}
)  \\
\end{array} 
\right.
\end{aligned}
\end{equation*}
for all $p'\in \mathcal{P}_1 \coloneq (0,1) \setminus \mathcal{P}_0= (0,a) \cup (b,1)$.
}
\newline
\textbf{Proof.} Now,
\begin{equation*}
\begin{aligned}
    \mathbb{P}
    \left[
    T_{m,r}>\epsilon
    \bigg| p'
    \right] &\geq
    \mathbb{P}
    \left[
    \tilde{T}_m>\epsilon+\sqrt{\frac{\mathrm{log}r}{r}}, 
    \left|
    T_{m,r}-\tilde{T}_m
    \right|<
    \sqrt{\frac{\mathrm{log}r}{r}} \bigg| p'
    \right] \\
    &\geq \mathbb{P}
    \left[
    \tilde{T}_m>
    \epsilon+
    \sqrt{\frac{\mathrm{log}r}{r}}
    \bigg| p'
    \right] +
    \mathbb{P}
    \left[
    \left|
T_{m,r}-\tilde{T}_m
    \right|
    <\sqrt{\frac{\mathrm{log}r}{r}}
    \bigg| p'
    \right]-1 \\
    &\geq 
    \mathbb{P}
    \left[
    \tilde{T}_m>
    \epsilon+
    \sqrt{\frac{\mathrm{log}r}{r}}
    \bigg| p'
    \right] +
    \left(
1-2\sum_{j=1}^m r^{-\frac{1}{2}}
    \right)-1 \\
&\geq \left\{
\begin{array}{ll}
    \left(
1- \frac{2m}{\sqrt{r}}
    \right)
      ,
       & p'\in \left(
       0,a-\epsilon-
       \sqrt{\frac{\mathrm{log}r}{r}}
       \right] \cup \left[
       b+\epsilon+
       \sqrt{\frac{\mathrm{log}r}{r}}
       ,1 
       \right) \\
       \left(
      \frac{b-p'-\epsilon-
      \sqrt{\frac{\mathrm{log}r}{r}}
      }{b-a}
      \right)^m
      - \frac{2m}{\sqrt{r}}
      ,
       & p' \in (a-\epsilon-
       \sqrt{\frac{\mathrm{log}r}{r}}
       ,a) \\
       \left(
      \frac{p'-\epsilon-
      \sqrt{\frac{\mathrm{log}r}{r}}
      -a}{b-a}
      \right)^m
      - \frac{2m}{\sqrt{r}}
      , & p' \in (b,b+\epsilon+
      \sqrt{\frac{\mathrm{log}r}{r}}
      )  \\
\end{array} 
\right. \\
&\geq 
\left\{
\begin{array}{ll}
    \left(
1- \frac{2m}{\sqrt{r}}
    \right)
      ,
       & p'\in \left(0,a-\epsilon-
       \sqrt{\frac{\mathrm{log}r}{r}}
       \right] \cup \left[
       b+\epsilon+
       \sqrt{\frac{\mathrm{log}r}{r}}
       ,1 
       \right) \\
       \left(
      \frac{b-a-\epsilon-
      \sqrt{\frac{\mathrm{log}r}{r}}
      }{b-a}
      \right)^m-
      \frac{2m}{\sqrt{r}}
      ,
       & p' \in (a-\epsilon-
       \sqrt{\frac{\mathrm{log}r}{r}}
       ,a) \\
       \left(
      \frac{b-\epsilon-
      \sqrt{\frac{\mathrm{log}r}{r}}
      -a}{b-a}
      \right)^m
      - \frac{2m}{\sqrt{r}}
      , & p' \in (b,b+\epsilon+
      \sqrt{\frac{\mathrm{log}r}{r}}
      )  \\
\end{array} 
\right. \\
&= 
\left\{
\begin{array}{ll}
    \left(
1- \frac{2m}{\sqrt{r}}
    \right)
      ,
       & p'\in \left(0,a-\epsilon-
       \sqrt{\frac{\mathrm{log}r}{r}}
       \right] \cup \left[
       b+\epsilon+
       \sqrt{\frac{\mathrm{log}r}{r}}
       ,1 
       \right) \\
       \left(1-
      \frac{\epsilon+
      \sqrt{\frac{\mathrm{log}r}{r}}
      }{b-a}
      \right)^m-
      \frac{2m}{\sqrt{r}}
      ,
       & p' \in (a-\epsilon-
       \sqrt{\frac{\mathrm{log}r}{r}}
       ,a) \\
       \left(1-
      \frac{\epsilon+
      \sqrt{\frac{\mathrm{log}r}{r}}
      }{b-a}
      \right)^m-
      \frac{2m}{\sqrt{r}}
      , & p' \in (b,b+\epsilon+
      \sqrt{\frac{\mathrm{log}r}{r}}
      )  \\
\end{array} 
\right. 
\end{aligned}
\end{equation*}
\newline
\newline
\newline
\textit{
\textbf{Theorem 4.} In the setting of \textit{Theorem \ref{Thm:lower_bound_power_realistic_h}}, under \textit{Assumptions \ref{Asm:null_uni_dist_h}} and \ref{Asm:alt_uni_dist_h}, when $\epsilon<\min \lbrace a,b-a,1-b \rbrace$,
    \begin{equation}
    \begin{aligned}
\mathbb{E}
\left[
\phi(p')
|p' \in \mathcal{P}_1
\right]
    = 
    \frac{2}{1-(b-a)}
    \left\lbrace
\left(
1-
\frac
{
\epsilon+
\sqrt{\frac{\mathrm{log}r}{r}}
}{b-a} 
\right)^m -
1
    \right\rbrace
    \left(
\epsilon+
\sqrt{\frac{\mathrm{log}r}{r}}
    \right) +
    \left(
1- \frac{2m}{\sqrt{r}}
    \right).
\end{aligned}
\end{equation}
}
\newline
\newline
\textbf{Proof.} 
First, observe that, if $p' \sim \mathrm{Unif}(\mathcal{P}_1)$, then the PDF will be given by
\begin{equation*}
    f_{'1}(p')=
    \left\{
\begin{array}{ll}
    \frac{1}{1-(b-a)}
    ,
       & p'\in (0,a) \cup (b,1) \\
       0, & \mathrm{o/w}  \\
\end{array} 
\right.
\end{equation*}
\newline
Note that
\begin{equation*}
\begin{aligned}
    &\mathbb{E}_{p' \sim \mathrm{Unif}(\mathcal{P}_1)}[\phi(p')] \\
    &= 
    \int_{-\infty}^{\infty
    } 
    \phi(p')f_{'1}(p')dp' \\
    &= \int_{0}^a 
    \phi(p')\frac{1}{1-(b-a)}dp'
    +
    \int_{b}^1 
    \phi(p')\frac{1}{1-(b-a)}
    dp' \\
    &= 
    \frac{2}{1-(b-a)}
    \left\lbrace
\left(
1-
\frac{
\epsilon+
\sqrt{\frac{\mathrm{log}r}{r}}
}{b-a}
\right)^m
-
\frac{2m}{\sqrt{r}}
    \right\rbrace
    \left(
\epsilon+
\sqrt{\frac{\mathrm{log}r}{r}}
    \right)  \\
    &\quad \quad +
    \frac{1}{1-(b-a)}
    \left(
1- \frac{2m}{\sqrt{r}}
    \right)
    \left\lbrace
1-(b-a)-
2
\left(
\epsilon+
\sqrt{\frac{\mathrm{log}r}{r}}
\right)
    \right\rbrace
    \\
    &=
    \frac{2}{1-(b-a)}
    \left\lbrace
\left(
1-
\frac
{
\epsilon+
\sqrt{\frac{\mathrm{log}r}{r}}
}{b-a} 
\right)^m -
1
    \right\rbrace
    \left(
\epsilon+
\sqrt{\frac{\mathrm{log}r}{r}}
    \right) +
    \left(
1- \frac{2m}{\sqrt{r}}
    \right)
\end{aligned}
\end{equation*}
when $\epsilon< \min \lbrace a,b-a,1-b \rbrace$ and $r$ is sufficiently large (using \textit{Theorem \ref{Thm:lower_bound_power_realistic_h}}). 
\newline
\newline
\newline
\textit{
\textbf{Corollary 1.} 
Consider the setting of \textit{Theorem \ref{Thm:lower_bound_power_realistic_h}}. As $\epsilon \to 0$, $m,r \to \infty$ such that $r=\omega(m^2)$ and 
$\sqrt{\frac{\mathrm{log}r}{r}} \leq \epsilon$, $\epsilon - \sqrt{\frac{\mathrm{log}r}{r}} \leq b-a$,
, we have an  asymptotically valid and consistent sequence of tests. 
}
\newline
\textbf{Proof.} Note that if $\epsilon-\sqrt{\frac{\mathrm{log}r}{r}} \leq (b-a)$, $\epsilon \geq \sqrt{\frac{\mathrm{log}r}{r}}$, 
$m \to \infty$, $r \to \infty$, then
\begin{equation*}
    \left(
1- 
\frac{\epsilon-
\sqrt{\frac{\mathrm{log}r}{r}}
}{b-a}
    \right)^m \to 0. 
\end{equation*}
Thus, if $\epsilon-\sqrt{\frac{\mathrm{log}r}{r}} \leq (b-a)$, $\epsilon \geq \sqrt{\frac{\mathrm{log}r}{r}}$, 
$m \to \infty$, $\frac{r^2}{m} \to \infty$,
\begin{equation*}
    \left(
1- 
\frac{\epsilon-
\sqrt{\frac{\mathrm{log}r}{r}}
}{b-a}
    \right)^m +
    \frac{2m}{\sqrt{r}}
    \to 0.
\end{equation*}
Using \textit{Theorem \ref{Thm:lower_bound_on_real_size_h}}, under the given conditions,
$\mathbb{P}[T_{m,r}>\epsilon|p'] \to 0$ for all $p' \in \mathcal{P}_0$. 
\newline
Also, note that, under the given conditions, $\phi(p') \to 1$ for all $p' \in \mathcal{P}_1$. 
\newline
\newline
\textit{
\textbf{Corollary 2.} 
As budget $\nu \to \infty$ such that $\tilde{m} \to \infty$ and $\tilde{r} \to \infty$, \textit{Algorithm \ref{Algo:TestOptimal}} yields an asymptotically valid sequence of tests. 
}
\newline
\textbf{Proof.} As $\tilde{m} \to \infty$, $\tilde{r} \to \infty$, 
\newline
$(\hat{b}-\hat{a}) \to^P (b-a)$ and hence 
\begin{equation*}
\left(
    1-\frac{1}{\hat{b}-\hat{a}}
    \left(
\epsilon-\sqrt{\frac{\mathrm{log}r}{r}}
    \right)
\right)^m + \frac{2m}{\sqrt{r}}
\to 
\left(
    1-\frac{1}{b-a}
    \left(
\epsilon-\sqrt{\frac{\mathrm{log}r}{r}}
    \right)
\right)^m + \frac{2m}{\sqrt{r}},
\end{equation*}
and thus the approximate validity constraint approaches the true validity constraint at level of significance $\alpha$. Thus, Algorithm \ref{Algo:TestOptimal} yields an asymptotically valid sequence of tests. 

\vskip 0.2in

\end{document}